\documentclass[a4paper,11pt,reqno,twoside]{amsart}
\usepackage{amsmath}
\usepackage{amsfonts}
\usepackage{amssymb}
\usepackage{amsthm}

\theoremstyle{plain}
\newtheorem{thm}{Theorem}[section]
\newtheorem{cor}[thm]{Corollary}
\newtheorem{lemma}[thm]{Lemma}
\newtheorem{prop}[thm]{Proposition}
\newtheorem{defn}[thm]{Definition}
\theoremstyle{definition}
\newtheorem{remark}[thm]{Remark}
{\theoremstyle{definition} 
                        \newtheorem*{example}{Example}
                        }

\numberwithin{equation}{section}

\newcommand{\N}{\mathbb{N}}                         
\newcommand{\C}{\mathbb{C}}                         
\newcommand{\Z}{\mathbb{Z}}                         

\newcommand{\h}{\mathfrak{h}}                       
\newcommand{\hs}{\mathfrak{h}^\ast}
\newcommand{\Mh}{M_{\h^*}}                          
\newcommand{\Dh}{{D_{\hs}}}
\newcommand{\A}{\mathcal{A}}                        
\newcommand{\B}{\mathcal{B}}                        
\newcommand{\F}{\mathcal{F}}
\newcommand{\U}{{\mathcal{U}}}                      
\newcommand{\opp}{{\mathrm{opp}}}                   
\newcommand{\cop}{{\mathrm{cop}}}                   
\newcommand{\lr}{{\mathrm{lr}}}                     
\newcommand{\ulambda}{\underline \lambda}
\newcommand{\umu}{\underline \mu}
\renewcommand{\emptyset}{\varnothing}               
\renewcommand{\det}{\mathrm{det}}                   
\newcommand{\1}{\mathbf{1}}                         
\newcommand{\id}{\mathrm{Id}}                       
\newcommand{\ihat}{\hat\imath}                      
\newcommand{\jhat}{\hat\jmath}
\newcommand{\totimes}{\widetilde \otimes}           
\newcommand{\da}{{\dagger}}
\newcommand{\sign}{\mathrm{sign}}                   

\newcounter{mylist}

\newenvironment{mylist}
{\begin{list} {(\roman{mylist}) }
   {\usecounter{mylist}
    \setlength{\labelsep}{0pt}
    \setlength{\leftmargin}{0pt}
    \setlength{\labelwidth}{0pt}
    \setlength{\topsep}{0.8ex plus 0.2ex minus 0ex}
    \setlength{\itemsep}{0ex plus 0.2ex minus 0ex}
    \setlength{\listparindent}{0pt}}}
{\end{list}}
\newenvironment{pair}[2]{ \langle#1, #2\rangle}                         

\newenvironment{qmindet}[2]{\xi^{\hat #1}_{\hat #2}}

\begin{document}
\date{August 4, 2005}
\title{The dynamical $U(n)$ quantum group}
\author{Erik Koelink and Yvette van Norden}
\thanks{The second author is supported by Netherlands Organisation for Scientific
Research (NWO) under project number 613.006.572.}
\address{Technische Universiteit Delft, EWI, DIAM, 
Postbus 5031, 2600 GA Delft, the Netherlands}
\email{h.t.koelink@ewi.tudelft.nl, y.vannorden@ewi.tudelft.nl}

\begin{abstract}
We study the dynamical analogue of the matrix algebra $M(n)$,
constructed from a dynamical $R$-matrix given by Etingof and Varchenko.
A left and a right corepresentation of this algebra, which can
be seen as analogues of the exterior algebra representation, are 
defined and this defines dynamical quantum minor determinants as
the matrix elements of these corepresentations. These elements
are studied in more detail, especially the action of the 
comultiplication and Laplace expansions. Using the Laplace
expansions we can prove that the dynamical quantum
determinant is almost central, and adjoining an inverse 
the antipode can be defined.
This results in the dynamical $GL(n)$
quantum group associated to the dynamical $R$-matrix. 
We study a $\ast$-structure leading to the dynamical $U(n)$ quantum
group, and we obtain results for the canonical pairing arising
from the $R$-matrix.
\end{abstract}

\maketitle


\section{Introduction}

Dynamical quantum groups have been introduced recently by 
Etingof and Varchenko \cite{EtinV2}, see 
the review paper by Etingof and Schiffmann \cite{EtingofSchiffmann} for
an overview and references to the literature, 
and related algebraic 
structures have been studied by Lu \cite{Lu}, Xu \cite{Xu}
in the context of deformations of Poisson groupoids, and
by Takeuchi \cite{Take}.  
Brzezi\`nski and Militaru \cite{BrzeM} compare the various
constructions of \cite{Lu}, \cite{Xu}, \cite{Take}. 
In this paper we stick to the definition of Etingof and
Varchenko \cite{EtinV2} with a slight modification as in 
\cite{KoelinkRosengren}. In order to keep the paper self-contained
as much as possible we recall the definition in 
section \ref{sec:dynMn}. We also recall the FRST-construction 
associated to a solution of the dynamical $R$-matrix,  
which gives a wealth of examples, and which we consider
explicitly for the trigonometric $R$-matrix 
in the $\mathfrak{gl}(n)$-case.

It is well-known that quantum groups have a natural 
link with special functions of basic hypergeometric type, 
and in \cite{KoelinkRosengren} it is shown that
this remains valid for the simplest example of 
a dynamical quantum group associated to the 
trigonometric $R$-matrix for $SL(2)$ and in \cite{KvNR},
see also \cite{KvN}, 
for the elliptic $R$-matrix for $SL(2)$ giving a 
dynamical quantum group theoretic interpretation
of elliptic hypergeometric series.  
In particular, \cite{KoelinkRosengren} gives a dynamical 
quantum group theoretic interpretation of 
Askey-Wilson and $q$-Racah polynomials having many
similarities to the interpretation of these 
polynomials on the (ordinary) quantum $SL(2)$ group
using the twisted primitive elements as introduced
by Koornwinder \cite{Koor}, see also \cite{NoumM}, 
\cite{KoelAAM}. This naturally suggests a link 
between these two approaches, and the link is established
by Stokman \cite{Stok} using the coboundary element
of Babelon, Bernard and Billey \cite{BabeBB}, a universal
element in the tensor product of the quantized 
universal algebra. This element also has a 
natural interpretation in the context of twisted
primitive elements as shown by Rosengren \cite{RoseCM}.
However, the coboundary element is only known for the
$\mathfrak{sl}(2)$-case, but there are conjectures 
about its form for the $\mathfrak{sl}(n)$ case, see
Buffenoir and Roche \cite{BuffR}. 
The notion of twisted primitive elements of Koornwinder \cite{Koor},
and especially its generalization to co-ideals, 
has turned out to be enormously fruitful for the interpretation
of special functions of one or many variables as spherical
functions on quantum groups or quantum symmetric spaces, see 
e.g. \cite{Dijk}, \cite{DijkN}, \cite{Letz}, \cite{Noum}, 
\cite{StokmanOb}. 

As one of the highlights of the application of 
Lie theory to special functions 
we mention the group theoretic derivation of the 
addition formula for Jacobi 
polynomials as obtained by Koornwinder \cite{KoorIM}
by working on the symmetric space $U(n)/U(n-1)$ 
and establishing the spherical and associated
spherical elements in terms of Jacobi polynomials,
see also Askey \cite[Lecture 4]{Aske} for a nice introduction.
In the case $n=2$ this gives the addition formula
for Legendre polynomials. For $q$-analogues of 
addition formulas for the Legendre polynomials 
see the overview \cite{KoelFI}. 
In the quantum group setting, Noumi, Yamada and Mimachi 
\cite{NoumiYM} established the little
$q$-Jacobi polynomials as spherical functions, and
Floris \cite{Flor} calculated the associated
spherical elements in terms of little $q$-Jacobi
polynomials and derived an addition formula. 
On the other hand, using the
notion of co-ideals,  
Dijkhuizen and Noumi \cite{DijkN} established 
Askey-Wilson polynomials as spherical functions
on a quantum analogue of $U(n)/U(n-1)$. 

In light of the above it is natural 
to ask for the spherical (and associated spherical) 
elements on the dynamical quantum
group analogue of $U(n)/U(n-1)$, and if 
a precise link to special functions can be established.
For this we need to study the dynamical 
$U(n)$ quantum group more closely, and in 
a previous paper \cite{KvN} we have 
studied general aspects of dynamical quantum groups
for this purpose. In case $n=2$ 
\cite{KoelinkRosengren} shows that the algebraic
approach to quantum groups as discussed by
Dijkhuizen and Koornwinder \cite{DijkK} is 
applicable, and we expect this to hold true for
general $n$.  
This paper serves as a first step in this specific
programme 
by defining the dynamical $U(n)$ quantum 
group and studying some of its 
elementary properties. In a future paper 
its corepresentation theory and (associated) 
spherical functions have to be studied. 

The general theory provides us with a dynamical
analogue of the algebra of functions on the 
space of $n\times n$-matrices. In section 
\ref{sec:dynMn} we recall the algebraic notions 
and FRST-construction 
of Etingof and Varchenko \cite{EtingofVarchenko}, \cite{EtinV2}, 
see also \cite{EtingofSchiffmann},  
and this gives an explicit presentation by
generators and relations for this dynamical
analogue. In order to make this into a dynamical
quantum group, we need to equip (a suitable
extension) of this algebra with an antipode.
For this purpose we study the dynamical analogues
of the minor determinants, which are introduced 
as matrix elements of corepresentations which
are analogues of the natural representation in 
the exterior algebra. There exist a left and
right corepresentation, and we show that the
matrix elements, i.e. the dynamical 
quantum minor determinants, are equal 
using an identity for Hall-Littlewood polynomials.
In particular, this gives a dynamical quantum
determinant. 
This is done in section  \ref{sec:extercorep}. 
In section \ref{sec:Laplace}, we continue the
study of these dynamical quantum minor
elements and we discuss the appropriate analogues
of the Laplace expansions. In section \ref{sec:dynGLn} 
we show how the Laplace
expansions imply that the dynamical
quantum determinant is almost central, and
localizing we find the dynamical
$GL(n)$ quantum group for which we give an explicit
expression for the antipode.  The treatment of
dynamical quantum minor elements, the Laplace
expansions, and the extension to a $\h$-Hopf algebroid
is very much motivated by the paper 
by Noumi, Yamada and Mimachi \cite{NoumiYM}. 
We also introduce 
a $\ast$-structure, so that we obtain 
the dynamical  $U(n)$ quantum group in section \ref{sec:dynUn}.
Finally, in section \ref{sec:pairing} we study the
natural pairing, as introduced by
Rosengren \cite{r}, see also \cite{KvN}, for the case of the
dynamical $GL(n)$ quantum group and the dynamical
$U(n)$ quantum group. 

\emph{Acknowledgement.} We thank Hjalmar
Rosengren and Jasper Stokman for useful
discussions.


\section{The dynamical analogue of the matrix algebra $M(n)$}\label{sec:dynMn}

In this section we give the general definitions of the theory of
dynamical quantum groups and we recall the generalized
FRST-construction. To define the $\h$-bialgebroid $\F_R(M(n))$ we
apply this construction to a solution of the quantum dynamical
Yang-Baxter equation (QDYBE).

Let $\h$ be a finite dimensional complex vector space, viewed as a
commutative Lie algebra, with dual space $\h^*$. Let
$V=\bigoplus_{\alpha\in\h^*} V_{\alpha}$ be a diagonalizable
$\h$-module. The quantum dynamical Yang-Baxter equation is given by
\begin{equation}\label{eq:QDYBE}
R^{12}(\lambda-h^{(3)})R^{13}(\lambda)R^{23}(\lambda-h^{(1)})=
R^{23}(\lambda)R^{13}(\lambda-h^{(2)})R^{12}(\lambda).
\end{equation}
Here $R:\h^*\to \mathrm{End}(V\otimes V)$ is a meromorphic function,
$h$ indicates the action of $\h$ and the upper indices are
leg-numbering notation for the tensor product. For instance,
$R^{12}(\lambda-h^{(3)})$ denotes the operator
$R^{12}(\lambda-h^{(3)})(u\otimes v \otimes w)=
R(\lambda-\mu)(u\otimes v) \otimes w$ for $w\in V_\mu$. An
$R$-matrix is a solution of the QDYBE \eqref{eq:QDYBE} which is
$\h$-invariant.

In the example we study, we identify $\h\cong\h^*\cong \C^n$ and
take $V$ an $n$-dimensional vector space with basis $\{v_1, \ldots,
v_n\}$. The $R$-matrix $R:\h^*\to \mathrm{End}_\h (V\otimes V)$ we
consider is given by
\begin{equation}\label{eq:RmatrixMn}
\begin{split}
R(\ulambda)= q\sum_{a=1}^n E_{aa}\otimes E_{aa}&+ \sum_{a<b}
E_{aa}\otimes E_{bb}  + \sum_{a>b} g(\lambda_{a}-\lambda_{b})
E_{aa}\otimes
 E_{bb}\\ &+ \sum_{a\neq b} h_0(\lambda_{a}-\lambda_{b}) E_{ba}\otimes E_{ab},
 \end{split}
\end{equation}
where $\ulambda= (\lambda_1, \ldots, \lambda_{n})$, $E_{ab}\in
\mathrm{End}(V)$ such that $E_{ab}v_c=\delta_{bc}v_a$ and the
meromorphic functions $h_0$ and $g$ are given by
\begin{equation}\label{eq:defnh0g}
h_0(\lambda)= \frac{q^{-1}-q}{q^{-2\lambda}-1}, \qquad
 g(\lambda)=\frac{(q^{-2\lambda}-q^{-2}) (q^{-2\lambda}-q^2)}
{(q^{-2\lambda}-1)^2}.
\end{equation}
Etingof and Varchenko \cite{EtinV2} obtain this $R$-matrix as the
exchange matrix for the vector representation of $GL(n)$.

\subsection{$\h$-Hopf algebroids and the generalized FRST-construction}

We recall the definition of $\h$-Hopf algebroids, the algebraic notion
for a dynamical quantum groups, and the generalized FRST-construction.

Let $\h$ be a finite dimensional complex vector space, with dual
space $\h^*$. Denote by $\Mh$ the field of meromorphic functions on $\h^*$.
For $\alpha\in\h^*$ we denote by $T_{\alpha}:\Mh\to\Mh$ the
automorphism $(T_\alpha f )(\lambda)=f(\lambda+\alpha)$ for all
$\lambda\in\h^*$.

\begin{defn}
A $\h$-algebra is a complex associative algebra $\A$ with $1$ which
is bigraded over $\h^*$, $\A=\bigoplus_{\alpha,\beta\in\h^*}
\A_{\alpha\beta}$, with two algebra embeddings $\mu_l$,
$\mu_r:\Mh\to \A_{00}$ (the left and right moment map) such that
$\mu_l(f)a=a\mu_l(T_{\alpha} f)$, $\mu_r(f)a=a\mu_r(T_{\beta}f)$,
for all $f\in\Mh$, $a\in\A_{\alpha\beta}$.

A morphism of $\h$-algebras is an algebra homomorphism which
preserves the bigrading and the moment maps.
\end{defn}

Let $\A$ and $\B$ be two $\h$-algebras. The matrix tensor product
$\A\tilde \otimes \B$ is the $\h^*$-bigraded vector space with $(\A
\tilde \otimes \B)_{\alpha\beta}= \bigoplus_{\gamma\in\h^*}
(\A_{\alpha\gamma}\otimes_{M_{\h^*}}\B_{\gamma\beta})$, where
$\otimes_{M_{\h^*}}$denotes the usual tensor product modulo the
relations
    \begin{equation}\label{eqn:tensor}
\mu_r^\A(f)a\otimes b= a\otimes \mu_l^\B(f)b, \;\mbox{for all } a\in
\A, b\in \B, f \in M_{\h^*}.
    \end{equation}
The multiplication $(a\otimes b)(c\otimes d)=ac\otimes bd$ for $a$,
$c\in \A$ and $b$, $d\in \B$ and the moment maps
$\mu_l(f)=\mu_l^\A(f)\otimes 1$ and $\mu_r(f)=1\otimes \mu_r^\B(f)$
make $\A\tilde \otimes \B$ into a $\h$-algebra.

\begin{example}
Let $\Dh$ be the algebra of difference operators acting on  $M_{\h^*}$,
consisting of the operators $\sum_i f_i T_{\beta_i}$, with $f_i\in
M_{\h^*}$ and $\beta_i\in\h^*$. This is a $\h$-algebra with the
bigrading defined by $fT_{-\beta} \in (\Dh)_{\beta\beta}$ and both
moment maps equal to the natural embedding.

For any $\h$-algebra $\A$, there are canonical isomorphisms $\A
\cong \A\tilde\otimes \Dh \cong \Dh\tilde\otimes \A$, defined by
    \begin{equation}\label{eqn:identDh}
x\cong x \otimes T_{-\beta} \cong T_{-\alpha}\otimes x, \; \mbox{for
all } x\in \A_{\alpha\beta}.
    \end{equation}
The algebra $\Dh$ plays the role of the unit object in the category
of $\h$-algebras.
\end{example}

\begin{defn}
A $\h$-bialgebroid is a $\h$-algebra $\A$ equipped with two
$\h$-algebra homomorphisms $\Delta:\A\to \A\tilde\otimes \A$ (the
comultiplication) and $\varepsilon:\A \to \Dh$ (the counit) such
that $(\Delta\otimes \id)\circ\Delta= (\id\otimes
\Delta)\circ\Delta$ and $(\varepsilon\otimes \id)\circ
\Delta=\id=(\id\otimes \varepsilon)\circ\Delta$ (under the
identifications \eqref{eqn:identDh}).
\end{defn}

For the definition of the antipode we follow
\cite{KoelinkRosengren}.

\begin{defn}\label{def:antipode}
A $\h$-Hopf algebroid is a $\h$-bialgebroid $\A$ equipped with a
$\C$-linear map $S:\A\to \A$, the antipode, such that
$S(\mu_r(f)a)=S(a)\mu_l(f)$ and $S(a\mu_l(f))=\mu_r(f)S(a)$ for all
$a\in \A$, $f\in M_{\h^*}$, and
    \begin{equation}\label{eq:defS}
    \begin{split}
& m\circ (\id\otimes S)\circ\Delta(a)= \mu_l(\varepsilon(a)1),\; \mbox{for all } a\in \A,\\
 &m\circ(S\otimes \id)\circ\Delta(a)=\mu_r(T_\alpha
 (\varepsilon(a)1)), \; \mbox{for all } a\in \A_{\alpha\beta},
    \end{split}
    \end{equation}
where $m:\A\times \A\to \A$ denotes the multiplication and
$\varepsilon(a)1$ is the result of applying the difference operator
$\varepsilon(a)$ to the constant function $1\in M_{\h^*}$.
\end{defn}

If there exists an antipode on a $\h$-bialgebroid, it is unique.
Furthermore, the antipode is anti-multiplicative,
anti-comultiplicative, unital, counital and interchanges the moment
maps $\mu_l$ and $\mu_r$, see \cite[Prop. 2.2]{KoelinkRosengren}.
In Definition \ref{def:antipode} the maps $m\circ (\id\otimes S)$
and $m\circ(S\otimes \id)$ are well-defined on $\A\tilde\otimes\A$,
see \cite{KvN}. 

\begin{example}\label{ex:oppcoop}
\begin{mylist}
\item We can equip $\Dh$ with a $\h$-Hopf algebroid structure with
comultiplication $\Delta:\Dh\to \Dh\tilde\otimes\Dh \cong\Dh$ the
canonical isomorphism, counit $\varepsilon:\Dh\to\Dh$ the identity
and antipode defined by $S(fT_{\alpha})= T_{-\alpha}\circ f$.
\item For a $\h$-Hopf algebroid $\A$ with invertible antipode, 
the opposite and co-opposite are
also $\h$-Hopf algebroids. The opposite algebra $\A^\opp$ is the
algebra $\A$ with opposite multiplication. Then we equip $\A^\opp$
with a $\h$-Hopf algebroid structure by defining
$(\A^\opp)_{\alpha\beta}=\A_{-\alpha,-\beta}$,
$\mu_l^\opp=\mu_l^\A$, $\mu_r^\opp=\mu_r^\A$,
$\Delta^\opp=\Delta^\A$, $\varepsilon^\opp=S^{\Dh}\circ
\varepsilon^\A$ and $S^\opp=(S^\A)^{-1}$. The co-opposite algebra
$\A^\cop$ has the same algebra structure but $\mu_l^\cop=\mu_r^\A$,
$\mu_r^\cop=\mu_l^\A$, $(\A^\cop)_{\alpha\beta}=\A_{\beta\alpha}$
and $\Delta^\cop=P\circ \Delta^\A$,
$\varepsilon^\cop=\varepsilon^\A$, $S^\cop=(S^\A)^{-1}$, where $P$
is the flip operator.
\end{mylist}
\end{example}

Let $\lambda\mapsto \overline \lambda$ be a  complex
conjugation on $\h^*$, and denote $\overline
f(\lambda)=\overline{f(\overline \lambda)}$ for all $f\in \Mh$.

\begin{defn}\label{def:hstarHopfalgebroid}
A $\h$-$\ast$-bialgebroid $\A$ is a $\h$-bialgebroid equipped 
with a $\ast$-operator, i.e. a $\C$-antilinear antimultiplicative
involution such that $\mu_l(f)=\mu_l(\bar f)$ and 
$\mu_r(f)=\mu_r(\bar f)$, such that $(\ast\otimes\ast)\circ\Delta =\Delta\circ\ast$
and $\varepsilon\circ\ast = \ast_{D_{\h^\ast}}\circ \varepsilon$, where 
$\ast=\ast_{D_{\h^\ast}}$ on $D_{\h^\ast}$ is defined by 
$(f T_{\alpha})^\ast= T_{-\bar\alpha}\circ \bar f$.
\end{defn}

A $\h$-Hopf $\ast$-algebroid is a $\h$-Hopf algebroid that is a $\h$-$\ast$-bialgebroid and has an invertible antipode. Then,
see \cite{KoelinkRosengren}, $S\circ\ast$ is an involution. 

Until this point we have seen only the example $D_{\h^\ast}$ of 
a $\h$-bialgebroid. 
The generalized FRST-construction provides many examples of
$\h$-bialgebroids from $R$-matrices, see \cite{EtingofVarchenko},
\cite{EtingofSchiffmann}, \cite{FelderVarchenko},
\cite{KoelinkRosengren}. We recall the construction and we
apply the construction to the $R$-matrix in \eqref{eq:RmatrixMn} to obtain
the main object of study for this paper. 

Let $\h$ and $M_{\h^*}$ be as before,
$V=\bigoplus_{\alpha\in\h^*}V_\alpha$ be a finite-dimensional
diagonalizable $\h$-module and $R:\h^* \to \mathrm{End}_\h(V\otimes
V)$ a meromorphic function that commutes with the $\h$-action on
$V\otimes V$. Let $\{v_x\}_{x\in X}$ be a homogeneous basis of $V$,
where $X$ is an index set. Write $R^{a b}_{x y}(\lambda)$ for the
matrix elements of $R$,
    \begin{eqnarray*}
R(\lambda)(v_a\otimes v_b)=\sum_{x,y\in X} R^{a b}_{x y}(\lambda)
v_x\otimes v_y,
    \end{eqnarray*}
and define $\omega:X\to \h^*$ by $v_x\in V_{\omega(x)}$. Let $\A_R$
be the unital complex associative algebra generated by the elements
$\{L_{xy}\}_{x,y\in X}$ together with two copies of $M_{\h^*}$,
embedded as subalgebras. The elements of these two copies will be
denoted by $f(\lambda)$ and $f(\mu)$, respectively. The defining
relations of $\A_R$ are $f(\lambda)g(\mu) = g(\mu)f(\lambda)$,
$f(\lambda)L_{x y} = L_{x y} f(\lambda+\omega(x))$ and $f(\mu)L_{x
y} =L_{x y} f(\mu+\omega(y))$ for all $f$, $g\in M_{\h^*}$, together
with the RLL-relations
    \begin{eqnarray}\label{eq:RLL}
\sum_{x,y\in X} R^{x y}_{a c}(\lambda)L_{x b}L_{y d}= \sum_{x,y\in
X} R_{x y}^{b d}(\mu)L_{c y}L_{a x},
    \end{eqnarray}
for all $a$, $b$, $c$, $d\in X$. The bigrading on $\A_R$ is defined
by $L_{x y}\in \A_{\omega(x),\omega(y)}$ and $f(\lambda)$,
$f(\mu)\in \A_{00}$. The moment maps defined by
$\mu_l(f)=f(\lambda)$, $\mu_r(f)=f(\mu)$ make $\A_R$ into a
$\h$-algebra. The $\h$-invariance of $R$ ensures that the bigrading
is compatible with the RLL-relations \eqref{eq:RLL}. Finally the
counit and comultiplication defined by
    \begin{gather}\label{eqn:counit}
\varepsilon(L_{a b})=\delta_{a b} T_{-\omega(a)},\quad
\varepsilon(f(\lambda))=\varepsilon(f(\mu))=f, \\
 \label{eqn:comult}
 \Delta(L_{a b})=\sum_{x\in X} L_{a x}\otimes L_{x b}, \quad
\Delta(f(\lambda))=f(\lambda)\otimes 1, \quad
\Delta(f(\mu))=1\otimes f(\mu),
    \end{gather}
equip $\A_R$ with the structure of a $\h$-bialgebroid, see
\cite{EtingofVarchenko}.

\subsection{The dynamical analogue of the algebra $M(n)$}

Now, we apply the generalized FRST-construction to the $R$-matrix
\eqref{eq:RmatrixMn} in order to define the $\h$-bialgebroid
$\F_R(M(n))$. Let $X=\{1, \ldots, n\}$ and define $\omega:X\to \h^*$
by $i\mapsto e_i$, where $e_i$ is the $i$-th unit vector of $\C^n$.
Let $h(\lambda)=q-h_0(\lambda)$ where $h_0$ is defined as in
\eqref{eq:defnh0g}, so
\begin{equation}\label{eq:defnh}
h(\lambda)= q\frac{(q^{-2\lambda}-q^{-2})}{(q^{-2\lambda}-1)}.
\end{equation}

\begin{defn}\label{def:dynMn}
The $\h$-algebra $\F_R(M(n))$ is the algebra generated by the
elements $t_{ij}$, $i$, $j\in\{1,2,\ldots,n\}$ together with two
copies of $\Mh$, denoted by $f(\ulambda)=f(\lambda_1, \ldots,
\lambda_n)$ and $f(\umu)=f(\mu_1,\ldots, \mu_n)$, embedded as
subalgebras. Then the defining relations are given by
$f_1(\ulambda)f_2(\umu) = f_2(\umu)f_1(\ulambda)$,
  \begin{equation}\label{eq:commdynMnftij}
    \begin{split}
  f(\ulambda)t_{ij} =t_{ij} f(\ulambda+\omega(i)),\quad
      f(\umu)t_{ij} =t_{ij} f(\umu+\omega(j)),
    \end{split}
    \end{equation}
with $f$, $f_1$, $f_2\in \Mh$, together with the RLL-relations
\begin{equation}\label{eq:RLLdynMn}
\begin{split}
  &h(\mu_{b}-\mu_{d})t_{ab}t_{ad}= t_{ad}t_{ab},\; \mbox{ for all } b<d,\\
  &h(\lambda_{c}-\lambda_{a})t_{cb}t_{ab}=t_{ab}t_{cb},\; \mbox{ for all } a<c,\\
  &t_{ab}t_{cd}= t_{cd}t_{ab} + [h(\lambda_{c}-\lambda_{a})-h(\mu_{b}-\mu_{d})]t_{cb}t_{ad},
         \;  \mbox{ for all }a<c, b<d\\
  & g(\mu_{b}-\mu_{d})t_{ab}t_{cd}=  g(\lambda_{a}-\lambda_{c})t_{cd}t_{ab}\\
    &\qquad \qquad\qquad\qquad + [h(\mu_{d}-\mu_{b})-h(\lambda_{a}-\lambda_{c})]t_{ad}t_{cb},
          \; \mbox{ for all } a<c, b<d,
\end{split}
\end{equation}
The bigrading $\F_R(M(n))=\bigoplus_{m,p\in \N^n} \F_{mp}$ is
defined on the generators by $f(\ulambda), f(\umu)\in\F_{00}$,
$t_{ij} \in\F_{\omega(i),\omega(j)}$ and the moment maps are given
by $\mu_l(f)=f(\ulambda)$, $\mu_r(f)=f(\umu)$. By defining the
comultiplication $\Delta: \F_R(M(n))\to \F_R(M(n))\tilde \otimes
\F_R(M(n))$ and counit $\varepsilon: \F_R(M(n))\to \Dh$ on the
generators by
 \begin{equation}
\Delta (t_{ij})= \sum_{k=1}^n t_{ik}\otimes t_{kj}, \qquad
 \Delta (f(\ulambda)) = f(\ulambda)\otimes 1,\qquad \Delta(f(\umu)) =1\otimes
 f(\umu),
   \end{equation}
and $\varepsilon(t_{ij})= \delta_{ij}T_{-\omega(i)}$,
$\varepsilon(f(\ulambda))=\varepsilon(f(\umu))=f$ and extended as
algebra homomorphisms we equip $\F_R(M(n))$ with the structure of a
$\h$-bialgebroid.
\end{defn}

\begin{remark}\label{rem:2inn}
The case $n=2$ and restricting to functions depending only
on $\lambda_1-\lambda_2$ gives back the case
studied in \cite{KoelinkRosengren}.
\end{remark}

As in \cite{NoumiYM} for the quantum case and
in \cite{KoelinkRosengren} for $n=2$, we can give a linear basis for
$\F_R(M(n))$. The proof is more involved since we use relations for
the functions $h$, $g$. Proposition \ref{prop:basisdynMn} is
stated for later reference. 

\begin{prop}\label{prop:basisdynMn}
For every $n\times n$-matrix $A$ we denote $t^A =
t_{11}^{a_{11}}t_{12}^{a_{12}}\cdots
t_{1n}^{a_{1n}}t_{21}^{a_{21}}\cdots t_{nn}^{a_{nn}}$. Then $\{t^A:
A\in M_n(\N)\}$ forms a basis over $\Mh\otimes \Mh$ for the vector
space $\F_R(M(n))$.
\end{prop}

\begin{proof}
This follows from the diamond lemma, see \cite{Bergman}. First we
introduce a total ordering $\prec$ by $t^A\prec t^B$ if $\sum_{i,j}
a_{ij}< \sum_{i,j} b_{ij}$ and in case $\sum_{i,j} a_{ij}=
\sum_{i,j} b_{ij}$ we use the lexicographical ordering on $(a_{11},
a_{12}, \ldots, a_{1n}, a_{21}, \ldots a_{2n}, a_{31}, \ldots,
a_{nn})$.

We have the following reduction system, which is compatible with the
introduced total order. Assume $i<j$, $k<l$,
\begin{equation*}
\begin{split}
     t_{il}t_{ik}\; \mapsto \;& h(\mu_k-\mu_l) t_{ik}t_{il},\\
     t_{jk}t_{ik} \;\mapsto \;& h(\lambda_j-\lambda_i)^{-1} t_{ik}
        t_{jk},\\
     t_{jl} t_{ik} \;\mapsto \;& [h(\lambda_j-\lambda_i)^{-1}+h(\mu_l-\mu_k)g(\lambda_i-\lambda_j)^{-1}] t_{il}t_{jk}
      +g(\mu_k-\mu_l)g(\lambda_i-\lambda_j)^{-1}t_{ik}t_{jl}\\
     t_{jk} t_{il} \;\mapsto \;& g(\lambda_i-\lambda_j)^{-1}t_{il}t_{jk}
      +  (h(\lambda_j-\lambda_i)^{-1}-h(\mu_k-\mu_l)g(\lambda_i-\lambda_j)^{-1})t_{ik}t_{jl}.\\
\end{split}
\end{equation*}
To simplify the coefficients on the right hand side we use
$h(\lambda)-h(-\mu)=h(\mu)-h(-\lambda)$ and
$g(\mu)-g(\lambda)=(h(\lambda)-h(-\mu))(h(\lambda)-h(\mu))$. If we
prove that the reduction system is resolvable, the lemma follows
from \cite[Thm 2.1]{Bergman}. There are 24 types of configuration to
be checked. The proof is straightforward using
\begin{equation}\label{eq:relationshg}
\begin{split}
h(-\lambda)= 1/ h(\lambda+1),& \quad
 g(-\lambda)= g(\lambda),  \quad g(\lambda) =
 h(\lambda)h(-\lambda),
\end{split}
\end{equation}
and the identities $h(\lambda)-h(\mu)=h(-\mu)-h(-\lambda)$ and
\begin{equation*}
\begin{split}
h(\lambda)h(\lambda-1)-h(\mu)&h(\nu)+h(\mu-\nu)h(\nu)-h(\mu-\nu)h(\lambda)
 + h(\nu-\mu)h(\mu)-h(\nu-\mu)h(\lambda)=0.
\end{split}
\end{equation*}
for all $\lambda$, $\mu$, $\nu\in\h^*$.
\end{proof}

\section{Exterior corepresentations and 
dynamical quantum minor determinants}\label{sec:extercorep}

We continue with the study of some elementary corepresentations of $\F_R(M(n))$ analogous to the action of $M(n)$ on the exterior algebra
of $\C^n$.  Using these corepresentations we find the dynamical
determinant in $\F_R(M(n))$. 
First we recall the general definition of a corepresentation of a
$\h$-bialgebroid on a $\h$-space, see \cite{KoelinkRosengren}. We
introduce the notion of $\h$-comodule algebras.

\begin{defn}
A $\h$-space is a vector space over $M_{\h^*}$ which is also a
diagonalizable $\h$-module, $V=\bigoplus_{\alpha\in\h^*}V_\alpha$,
with $M_{\h^*}V_\alpha\subseteq V_\alpha$ for all $\alpha\in\h^*$. A
morphism of $\h$-spaces is a $\h$-invariant (i.e.~ grade preserving)
$M_{\h^*}$-linear map.
\end{defn}

In case we want to emphasize the dependence on $V$ we 
also write $fv=\mu_V(f)v$. 

We next define the tensor product of a $\h$-bialgebroid $\A$ and a
$\h$-space $V$. Put $V\tilde \otimes
\A=\bigoplus_{\alpha,\beta\in\h^*}
(V_{\alpha}\otimes_{\Mh}\A_{\alpha\beta})$ where $\otimes_{\Mh}$
denotes the usual tensor product modulo the relations $v\otimes
\mu_l(f)a=fv\otimes a$. The grading $V_{\alpha}\otimes_{\Mh}
\A_{\alpha\beta} \subseteq (V\otimes \A)_{\beta}$ for all $\alpha$
and  $f(v\otimes a)=v\otimes \mu_r(f)a$ make
$V\tilde\otimes \A$ into a $\h$-space. Analogously  $\A\tilde\otimes
V=\bigoplus_{\alpha,\beta\in\h^*}(\A_{\alpha\beta}
\otimes_{M_{\h^*}}V_\beta)$
where $\otimes_{M_{\h^*}}$ denotes the usual tensor product modulo
the relations $\mu_r(f) a\otimes v=a\otimes f v$. The grading
$\A_{\alpha\beta}\otimes_{M_{\h^*}}V_\beta \subseteq
(\A\tilde\otimes V)_\alpha$ and $f(a\otimes
v)=\mu_l(f)a\otimes v$, $a\in \A$, $v\in V$, $f\in M_{\h^*}$, make
$\A\tilde\otimes V$ into a $\h$-space.

\begin{defn}
A right corepresentation of a $\h$-bialgebroid $\A$ on a $\h$-space
$V$ is a $\h$-space morphism $\rho:V\to V\tilde\otimes \A$ such that
$(\id\otimes \Delta)\circ \rho=(\rho\otimes\id)\circ\rho$,
$(\id\otimes \varepsilon)\circ \rho=\id$. The first equality is in
the sense of the natural isomorphism $(V\tilde\otimes
\A)\tilde\otimes \A\cong V\tilde\otimes (\A\tilde\otimes \A)$ and in
the second identity we use the identification $V\cong V
\tilde\otimes \Dh$ defined by $v\otimes fT_{-\alpha}\cong f v$,
$f\in M_{\h^*}$, for all $v\in V_\alpha$.

A left corepresentation of a $\h$-bialgebroid $\A$ on a $\h$-space
$V$ is a $\h$-space morphism $\rho:V\to \A \tilde\otimes V$ such
that $(\Delta\otimes \id)\circ \rho=(\id\otimes\rho)\circ\rho$,
$(\varepsilon\otimes \id)\circ \rho=\id$.
\end{defn}

\begin{defn}\label{def:hcomodule}
Let $\A$ be a $\h$-bialgebroid and $V$ a $\h$-space. Then $V$ is a
right (left) \emph{$\h$-comodule algebra} for $\A$ if there exists a
right (left) corepresentation $R\colon V \to V\widetilde\otimes \A$
($L: V\to \A\totimes V$) such that
\begin{mylist}
\item \label{item:hcomodule} $V$ is an associative algebra such that
$\mu_V(f) v w = v \mu_V(T_\alpha f)w$ for $v\in V_\alpha$, $w\in V$, and
$V_\alpha V_\beta\subset V_{\alpha+\beta}$,
\item $R$ ($L$) is an algebra homomorphism.
\end{mylist}
If, moreover, $V$ is a unital algebra, we require $R$ ($L$) to be
unital.
\end{defn}

\begin{remark}\label{rem:hcomodule}
The algebra structure of $V\totimes \A$ is given by $(v\otimes
a)(w\otimes b)=vw\otimes ab$ for $v$, $w\in V$ and $a$, $b\in \A$.
For $v\in V_\alpha$, $w\in V_{\gamma}$ and $a\in \A_{\alpha\beta}$,
$b\in\A_{\gamma\delta}$, we have $(v\otimes a)(w\otimes b)=vw\otimes
ab\in V_{\alpha+\gamma}\otimes \A_{\alpha+\gamma, \beta+\delta}$
using \ref{item:hcomodule} which implies $(V\totimes \A)_\beta
(V\totimes \A)_\delta \subseteq (V\totimes \A)_{\beta+\delta}$. For
$v\in V_\alpha$,
\begin{equation*}
\begin{split}
\mu_{V\totimes A}(f) R(v)=(1\otimes \mu_r(f))R(v)  = R(v)(1\otimes
\mu_r(T_\alpha f))=R(v)\mu_{V\totimes A}(T_\alpha f).
\end{split}
\end{equation*}
So $R$ preserves the relation in \ref{item:hcomodule}. Recall that
by the $\Mh$-linearity of a corepresentation we have
$R(\mu_V(f)v)=\mu_{V\totimes A}(f) R(v)=(1\otimes \mu_r(f))R(v)$.
\end{remark}

Now we define the $\h$-space $W$ on which we construct a right
corepresentation of $\F_R(M(n))$. $W$ can be seen as the dynamical
analogue of the exterior algebra representation.

\begin{defn}
Let $W$ be the unital associative algebra generated by the elements
$w_i$, $i\in\{1,2,\ldots,n\}$ and a copy of $M_{\h^\ast}$ embedded
as a subalgebra, its elements denoted by $f(\underline{\lambda})$,
subject to the relations
\begin{equation}
\begin{split}\label{eq:relationsW}
w_i^2=0 \mbox{ for all } i,\quad w_j w_i=
-h(\lambda_{j}-\lambda_{i})w_i w_j \mbox{ for all }i<j,
\end{split}
\end{equation}
with $h$ defined by \eqref{eq:defnh} and $f({\ulambda}) w_i = w_i
f({\ulambda}+\omega(i))$ for all $f\in \Mh$.
\end{defn}

For an ordered subset  $I=\{i_1,\ldots,i_r\}$, $1\leq i_1< \ldots <i_r\leq n$, 
of $\{1,\ldots,n\}$ 
we use the convention  $w_I=w_{i_1}\cdots
w_{i_r}$, unless mentioned otherwise. Moreover, $\emptyset$ is 
an ordered subset and $w_{\emptyset}=1$ corresponding to
the case $r=0$.  The following lemma is easily proved.

\begin{lemma}\label{lem:Wbasis}
$\dim_{\Mh} W=2^n$ and a basis for $W$ is given by $\{w_I: I=\{i_1,
\ldots, i_r\}, i_1< \ldots < i_r, r=1,\ldots, n\}$. $W$ is a
$\h$-space with $\mu_W(f)=f(\ulambda)$ and $w_I\in W_{\omega(I)}$
for $w_I$ a basis element with $\omega(I)=\sum_{j=1}^r \omega(i_j)$.
\end{lemma}

Define $W^r=\mathrm{span}_{\Mh}\{w_I: \#I=r\}$. Then $W=\bigoplus_{r=0}^n W^r$ and $W^r W^s
\subset W^{r+s}$, with the convention that $W^r=\{0\}$ if $r>n$.

\begin{prop}\label{prop:extensionR}
Define $R(1)=1\otimes 1$, $R(w_i)= \sum_{j=1}^n w_j\otimes t_{ji}$.
Then $R$ extends uniquely to $R\colon W\to W \widetilde \otimes
\F_R(M(n))$ such that $W$ is a right $\h$-comodule algebra for
$\F_R(M(n))$.
\end{prop}

\begin{proof}
It is clear that $W$ satisfies the conditions of Definition
\ref{def:hcomodule}\ref{item:hcomodule}. To see that $R$ can be
extended uniquely to an algebra homomorphism we need to verify
\begin{equation*}
R(w_i) R(w_i)=0 \mbox{ for all } i,\quad R(w_j) R(w_i)=
-R(h(\lambda_{j}-\lambda_{i})) R(w_i) R(w_j) \mbox{ for all }i<j,
\end{equation*}
and $R(f(\ulambda)) R(w_i) = R(w_i) R(f(\ulambda+\omega(i)))$ for
all $f\in \Mh$. By definition of $R$ and the defining relations
\eqref{eq:relationsW} of $W$ we get
\begin{subequations}
\begin{align}
R(w_i)R(w_i)&=\sum_{j ,k} w_jw_k\otimes t_{ji}t_{ki}
\label{eq:proofR1} = \sum_{k>j} [w_jw_k\otimes t_{ji}t_{ki}
-h(\lambda_k-\lambda_j)w_jw_k \otimes t_{ki}t_{ji}]\\
 & \label{eq:proofR2} =\sum_{k>j} w_jw_k\otimes
 (t_{ji}t_{ki}-h(\lambda_k-\lambda_j)t_{ki}t_{ji})=0,
\end{align}
\end{subequations}
where we use the second relation of \eqref{eq:RLLdynMn} in the last
equality. Let us emphasize that the function $h$ should be
interpreted in \eqref{eq:proofR1} as $\mu_W(\ulambda\mapsto
h(\lambda_k-\lambda_j))$ and in \eqref{eq:proofR2} as
$\mu_l(\ulambda\mapsto h(\lambda_k-\lambda_j))$. Similarly we obtain
that the relation $ R(w_j) R(w_i)= -R(h(\lambda_{j}-\lambda_{i}))
R(w_i) R(w_j)$ for $i<j$ is equivalent to
\begin{equation*}
\sum_{l>k} w_kw_l\otimes (t_{kj}t_{li}
-h(\lambda_l-\lambda_k)t_{lj}t_{ki} +h(\mu_j-\mu_i)t_{ki}t_{lj} -
h(\mu_j-\mu_i)h(\lambda_l-\lambda_k)t_{li}t_{kj})=0.
\end{equation*}
Using Lemma \ref{lem:Wbasis} it remains to prove that
\begin{equation}\label{eq:relationproofR}
t_{kj}t_{li} -h(\lambda_l-\lambda_k)t_{lj}t_{ki}
+h(\mu_j-\mu_i)t_{ki}t_{lj} -
h(\mu_j-\mu_i)h(\lambda_l-\lambda_k)t_{li}t_{kj}=0,
\end{equation}
for $i<j$, $k<l$. To show this we multiply this equation by
$h(\lambda_l-\lambda_k)-h(\mu_i-\mu_j)$ and eliminate the products
$t_{kj}t_{li}$ and $t_{li}t_{kj}$ using the third and fourth
relation in \eqref{eq:RLLdynMn} respectively. Using
$h(\lambda)-h(-\mu)=h(\mu)-h(-\lambda)$ for all $\lambda, \mu$ we
obtain that the relation \eqref{eq:relationproofR} holds 
by \eqref{eq:relationshg}. Using the definition of $R$
and Remark \ref{rem:hcomodule} the last relation follows
analogously.

By the definition of the comultiplication and the counit on the
generators of $\F_R(M(n))$ of Definition \ref{def:dynMn} it
immediately follows that $(\id\otimes \Delta)\circ R(w_i)=(R\otimes
\id)\circ R(w_i)$ and $(\id \otimes\varepsilon)\circ R(w_i)=w_i$.
Since $R$, $\Delta$ and $\varepsilon$ are $\h$-algebra
homomorphisms, so are $(\id\otimes \Delta)\circ R$, $(R\otimes
\id)\circ R$ and $(\id\otimes \varepsilon)\circ R$. So the
equalities $(\id\otimes \Delta)\circ R=(R\otimes \id)\circ R$ and
$(\id\otimes \varepsilon)\circ R=\id$ hold on the generators and
hence on all of $W$. $R$ is a corepresentation of $\F_R(M(n))$ on
$W$ and $W$ is a $\h$-comodule algebra for $\F_R(M(n))$.
\end{proof}

For $I$ and $J$ ordered subsets with $\#I=\#J$ we 
define the elements $\xi^I_J$ as the corresponding
matrix elements; 
\begin{equation*}
R(w_J)=\sum_{\#I=\#J} w_I\otimes \xi^I_J.
\end{equation*}
We use the convention that $\xi^I_J=0$ for all $I$, $J$ such
that $\#I\neq \#J$. 
In the remainder of the paper use the convention that a
summation over subsets such as $\sum_{\# I=r}$ is a summation over
all ordered subsets $I$ such that $\# I=r$.

\begin{cor}\label{cor:comultxi} 
\begin{mylist}
\item $\displaystyle{\Delta(\xi^I_J)=\sum_{\#K=\#I} \xi^I_K\otimes \xi^K_J}$ and
$\varepsilon(\xi^I_J)=\delta_{IJ}T_{-\omega(I)}$ for all $I$, $J$ with $\# I=\# J$,
\item $\displaystyle{R(W^r)\subset W^r\totimes \F_R(M(n))}$.
\end{mylist}
\end{cor}

We call the matrix elements $\xi_J^I$ the {\it dynamical quantum
minor determinants} of $\F_R(M(n))$ with respect to the subsets $I$
and $J$. The element $\xi_{\{1,\ldots,n\}}^{\{1,\ldots,n\}}$ is
called the \emph{determinant} of $\F_R(M(n))$, and is also
denoted by $\det$.

This right corepresentation has a left analogue, 
a left $\h$-comodule algebra $V$ for
$\F_R(M(n))$. The proofs are analogous to the ones for the right
$\h$-comodule algebra $W$, and are skipped.

\begin{defn}
Let $V$ be the unital associative algebra generated by the elements
$v_i$, $i\in \{1,\ldots,n\}$ and a copy of $\Mh$ embedded as a
subalgebra, its elements denoted by $f(\ulambda)$, subject to the
relations
\begin{equation}
v_i^2=0, \mbox{ for all i}, \qquad v_i v_j=
-h(\lambda_{i}-\lambda_{j})v_j v_i, \mbox{ for all $i<j$}.
\end{equation}
and $f(\ulambda)v_i=v_i f(\ulambda+\omega(i))$ for all $f\in\Mh$.
\end{defn}

For an ordered subset $I=\{i_1, \ldots i_r\}$ with $1\leq i_1<\ldots< i_r\leq
n$ we denote by $v_I$ the ordered element $v_I=v_{i_r}\cdots
v_{i_1}\in V$. Let us emphasize that an element $v_I\in V$ has
reversed order compared to $w_I\in W$ by notational convention.

\begin{lemma}
$\dim_{\Mh} V=2^n$ and a basis for $V$ is given by $\{v_I: I=\{i_1,
\ldots ,i_r\}, i_1< \ldots <i_r, r=1,\ldots, n\}$. $V$ is a
$\h$-space with $\mu_V(f)=f(\ulambda)$ and $v_I\in V_{\omega(I)}$
for $v_I$ a basis element.
\end{lemma}

Define $V^r=\mathrm{span}_{\Mh}\{v_I: \#I=r\}$. Then $V= \bigoplus_{r=0}^n V^r$ and 
$V^rV^s\subset V^{r+s}$, with the convention that $V^r=\{0\}$ if $r>n$.

\begin{prop}\label{prop:extensionL}
Define $L(1)=1\otimes 1$, $L(v_i)=\sum_{j=1}^n t_{ij}\otimes v_j$.
Then $L$ extends uniquely to $L:V\to \F_R(M(n))\totimes V$ such that
$V$ is a left $\h$-comodule algebra for $\F_R(M(n))$.
\end{prop}

For ordered subsets $I$, $J$ with $\#I=\#J$ 
we define the elements $\eta^I_J$ by
\[L(v_I)=\sum_{\# J=\#I} \eta^I_J \otimes v_J\]
and $\eta^I_J=0$ for $\#I\not=\#J$. 
We denote the corresponding determinant by $\widetilde
\det=\eta_{\{1,\ldots,n\}}^{\{1,\ldots,n\}}$.

\begin{cor}
\begin{mylist}
\item $\displaystyle{\Delta(\eta^I_J)=\sum_{\#K=\#I} 
\eta^I_K\otimes \eta^K_J}$ and
$\varepsilon(\eta^I_J)=\delta_{IJ}T_{-\omega(I)}$ for  $I$, $J$ with $\#I=\#J$,
\item $L(V^r)\subset \F_R(M(n)) \totimes V^r$.
\end{mylist}
\end{cor}

We call the matrix elements $\eta^I_J$ the dynamical quantum minor
determinants of $\F_R(M(n))$ with respect to the subsets $I$ and
$J$. In Proposition \ref{thm:xi=eta} we prove that the dynamical
quantum minor determinants related to the right and left
corepresentation are equal, so we can speak of \emph{the} dynamical
quantum minor determinants of $\F_R(M(n))$, without mentioning right
or left. First we compute an explicit expression of the dynamical
quantum minor determinants which we use in the proof.

For any permutation $\sigma\in S_r$, $1\leq r\leq n$, and any 
ordered subset
$I=\{i_1,\ldots,i_r\}$, we
define the generalized sign function $S(\sigma,I)\in M_{\h^\ast}$ by
\begin{equation}
\begin{split}
S(\sigma,I)(\ulambda) &= \prod_{\{k<l: \sigma(k)>\sigma(l)\}}
-h(\lambda_{i_{\sigma(k)}}-\lambda_{i_{\sigma(l)}})  =
(-q)^{l(\sigma)} \prod_{\{k<l: \sigma(k)>\sigma(l)\}}
\frac{q^{-2\lambda_{i_{\sigma(k)}}} -
q^{-2}q^{-2\lambda_{i_{\sigma(l)}}}} {q^{-2\lambda_{i_{\sigma(k)}}}
- q^{-2\lambda_{i_{\sigma(l)}}}},
\end{split}
\end{equation}
where $l(\sigma)$ denotes the length of the permutation,  
$l(\sigma)= \# \{k<l: \sigma(k)>\sigma(l)\}$.

\begin{lemma}\label{lem:orderW}
For any permutation $\sigma\in S_r$ we have the following relation
in $W$;
\begin{equation*}
w_{i_{\sigma(1)}} \ldots w_{i_{\sigma(r)}} =\mu_W(S(\sigma,I))\,
w_I,
\end{equation*}
where $I=\{i_1,\ldots,i_r\}$ is ordered.
\end{lemma}

\begin{proof}
We prove by induction on $r$, for $r=2$ and $\sigma=\id$ it is
trivial. If $\sigma=(12)$ it is just \eqref{eq:relationsW} for
$j=i_{\sigma(1)}$, $i=i_{\sigma(2)}$. Denote by $I'$ the ordered
subset of $I$ defined by $I\setminus \{i_{\sigma(1)}\}$, then
\begin{equation*}
\begin{split}
w_{i_{\sigma(1)}}\cdots &w_{i_{\sigma(r+1)}} = w_{i_{\sigma(1)}}
\prod_{\substack{2\leq k<l\leq r+1\\ \sigma(k)>\sigma(l)}}
-h(\lambda_{i_{\sigma(k)}}-\lambda_{i_{\sigma(l)}}) w_{I'}\\
 & = \prod_{\substack{2\leq k<l\leq r+1\\ \sigma(k)>\sigma(l)}}
-h(\lambda_{i_{\sigma(k)}}-\lambda_{i_{\sigma(l)}})
\prod_{\substack{2\leq l\leq r+1\\ \sigma(1)>\sigma(l)}}
-h(\lambda_{i_{\sigma(1)}}-\lambda_{i_{\sigma(l)}}) w_I  =
\mu_W(S(\sigma,I))w_I,
\end{split}
\end{equation*}
since $w_{i_{\sigma(1)}}$ commutes with all functions in $\Mh$ which
are independent of $\lambda_{i_{\sigma(1)}}$.
\end{proof}

Using Lemma \ref{lem:orderW} we calculate the action of the
corepresentation $R$ on $w_{j_1}\ldots w_{j_r}$ 
for an arbitrary  unordered set $\{j_1,\ldots, j_r\}$. Then
\begin{equation}\label{eq:expansionR}
R(w_{j_1}\ldots w_{j_r}) = R(w_{j_1})\ldots R(w_{j_r}) =
\sum_{k_1=1}^n \ldots \sum_{k_r=1}^n w_{k_1}\ldots w_{k_r} \otimes
t_{k_1j_1}\ldots t_{k_rj_r},
\end{equation}
and there is only a non-zero contribution in the right hand side of
\eqref{eq:expansionR} if all $k_i\not= k_j$ for $i\not= j$. Let
$I=\{i_1,\ldots, i_r\}$ be ordered, then we see that the contribution on the
right hand side  of \eqref{eq:expansionR} containing the basis
element $w_I$ in the first leg of the tensor product is given for
those terms for which $\{i_1,\ldots, i_r\}=\{k_1,\ldots, k_r\}$ as
unordered sets. So there exists for each non-zero term in
\eqref{eq:expansionR} contributing to the term containing $w_I$ in
the first leg of the tensor product precisely one permutation
$\sigma\in S_r$ such that $k_p=i_{\sigma(p)}$. So the term
containing $w_I$ in the first leg of the tensor product equals
\begin{equation*}
\begin{split}
\sum_{\sigma\in S_r} w_{i_{\sigma(1)}}\ldots w_{i_{\sigma(r)}}
\otimes t_{i_{\sigma(1)}j_1}\ldots t_{i_{\sigma(r)}j_r}  & =
\sum_{\sigma\in S_r}\mu_W(S(\sigma,I))\ w_I
\otimes t_{i_{\sigma(1)}j_1}\ldots t_{i_{\sigma(r)}j_r} \\
& =  \sum_{\sigma\in S_r} w_I \otimes \mu_l(S(\sigma,I))\
t_{i_{\sigma(1)}j_1}\ldots t_{i_{\sigma(r)}j_r},
\end{split}
\end{equation*}
by Lemma \ref{lem:orderW} and Remark \ref{rem:hcomodule}.

\begin{prop}\label{prop:matrixelementsR}
Let $J$ be ordered with $r=\#J$, then $R (w_J) =
\sum_{\# I = \#J } w_I \otimes \xi^I_J$ with the dynamical quantum
minor determinants given by
\begin{equation}\label{eq:xi}
\xi^I_J =\mu_r (S(\rho,J)^{-1}) \sum_{\sigma\in S_r}
\mu_l(S(\sigma,I))\ t_{i_{\sigma(1)}j_{\rho(1)}}\ldots
t_{i_{\sigma(r)}j_{\rho(r)}},
\end{equation}
for any $\rho\in S_r$.
\end{prop}

\begin{proof}
By Lemma \ref{lem:orderW} and the discussion preceding this
proposition we obtain
\begin{equation*}
\begin{split}
\sum_{\#I=\#J} w_I\otimes \xi^I_J &= R(w_J) = (1\otimes
\mu_r(S(\rho,J)^{-1})) R(w_{j_{\rho(1)}}\cdots w_{j_{\rho(r)}})\\
 &=  (1\otimes \mu_r(S(\rho,J)^{-1}))\sum_{\#I=\#J} 
 w_I\otimes \mu_l(S(\sigma,I))
t_{i_{\sigma(1)}j_{\rho(1)}}\ldots t_{i_{\sigma(r)}j_{\rho(r)}}.
\end{split}
\end{equation*}
So, the proposition follows from Lemma \ref{lem:Wbasis}.
\end{proof}

\begin{cor}
Put $S(\sigma)=S(\sigma,\{1,\ldots,n\})$ for 
$\sigma\in S_n$, then for any  $\rho\in S_n$,
\begin{equation*}
\det = \mu_r (S(\rho)^{-1}) \sum_{\sigma\in S_n} \mu_l(S(\sigma))\
t_{\sigma(1)\rho(1)}\ldots t_{\sigma(n)\rho(n)}
\end{equation*}
\end{cor}

Analogously we obtain an explicit formula for the matrix elements
$\eta^I_J$ of $L$. We need to define another generalized sign
function $\tilde S$ depending on an ordered subset $I$, $\#I=r$,
and  a permutation $\sigma\in S_r$; 
\begin{equation*}
\tilde S(\sigma,I)(\ulambda) :=\prod_{\{k<l: \sigma(k)>\sigma(l)\}}
-h(\lambda_{i_{\sigma(l)}}-\lambda_{i_{\sigma(k)}})=
\frac{1}{S(\sigma,I)(\lambda+1)}, 
\end{equation*}
where we use $h(-\lambda)=1/h(\lambda+1)$ for the last equality.
Analogous to Lemma \ref{lem:orderW}, we have for any permutation
$\sigma \in S_r$ the following relation in $V$
\begin{equation*}
v_{i_{\sigma(r)}} \cdots v_{i_{\sigma(1)}} = \mu_V(\tilde
S(\sigma,I))v_I,
\end{equation*}
where $I=\{i_1,\ldots,i_r\}$ is an ordered subset and
$v_I=v_{i_r}\cdots v_{i_1}$. We get the analogous
statement of Proposition \ref{prop:matrixelementsR}.

\begin{prop}
Let $I=\{i_1,\ldots,i_r\}$ be an ordered subset, then
$L(v_I)=\sum_{\# J=\#I} \eta^I_J\otimes v_J$ with the dynamical
quantum minor determinants given by, for any $\rho\in S_r$,
\begin{equation}\label{eq:eta}
\eta^I_J =\mu_l(\tilde S(\rho,I)^{-1})\sum_{\sigma\in S_r}
\mu_r(\tilde S(\sigma,J))t_{i_{\rho(r)}j_{\sigma(r)}}\ldots
t_{i_{\rho(1)}j_{\sigma(1)}}.
\end{equation}
\end{prop}

We now relate the two sets of dynamical quantum minor
determinants. For this we need the following identity;
\begin{equation}
\label{eq:HLidentity} \sum_{\sigma\in S_r} \prod_{i<j}
\frac{x_{\sigma(i)}-tx_{\sigma(j)}} {x_{\sigma(i)}-x_{\sigma(j)}} =
\prod_{i=1}^r \frac{1-t^i}{1-t}
\end{equation}
for $r$ indeterminates $x_1,\ldots x_r$. This identity can be found
in Macdonald \cite[III.1, (1.4)]{Macdonald} as the identity
expressing that the Hall-Littlewood polynomials for the zero
partition gives $1$.

\begin{thm}\label{thm:xi=eta} $\xi^I_J=\eta^I_J$ in $\F_R(M(n))$.
\end{thm}

\begin{proof} The proof is based on the expressions 
\eqref{eq:xi} and \eqref{eq:eta}, which give the possibility
to write a suitable multiple of $\xi^I_J$ as a double
sum over $S_r$, which, by interchanging summations,
gives a multiple of $\eta^I_J$. The multiples turn
out to be equal. The details are as follows.

First we rewrite $\eta^I_J$. Define the longest element $\sigma_0\in
S_r$ by $\sigma_0=\bigl(\begin{smallmatrix}1 & 2 &\ldots & r\\ r&
r-1&\ldots & 1\end{smallmatrix}\bigr)$. By substituting $\rho\mapsto
\rho\sigma_0$ and $\sigma\mapsto \sigma\sigma_0$ in \eqref{eq:eta}
we get
\begin{equation*}
\eta^I_J=\prod_{\substack{m<p\\ \rho(m)< \rho(p)}}
-h(\lambda_{i_{\rho(m)}}-\lambda_{i_{\rho(p)}})^{-1}
\sum_{\sigma\in S_r} \prod_{\substack{k<l\\
\sigma(k)<\sigma(l)}}-h(\mu_{j_{\sigma(k)}}-\mu_{j_{\sigma(l)}})
t_{i_{\rho(1)}j_{\sigma(1)}}\cdots t_{i_{\rho(r)}j_{\sigma(r)}},
\end{equation*}
for any $\rho\in S_r$. Using this expression for $\eta^I_J$ and
\eqref{eq:xi} we compute
\begin{equation*}
\begin{split}
&\left( \sum_{\rho\in S_r}\prod_{k<l}
-h(\mu_{j_{\rho(k)}}-\mu_{j_{\rho(l)}}) \right) \xi^I_J\\
 &\quad =\sum_{\rho\in S_r}\sum_{\sigma\in S_r} \prod_{\substack{k<l\\
\rho(k)<\rho(l)}}-h(\mu_{j_{\rho(k)}}- \mu_{j_{\rho(l)}})
\prod_{\substack{k<l\\\sigma(k)>\sigma(l)}}-h(\lambda_{i_{\sigma(k)}}-\lambda_{i_{\sigma(l)}})
t_{i_{\sigma(1)}j_{\rho(1)}}\cdots t_{i_{\sigma(r)}j_{\rho(r)}} \\
 &\quad =  \sum_{\sigma\in S_r}\sum_{\rho\in S_r}
 \prod_{\substack{k<l}}-h(\lambda_{i_{\sigma(k)}}-\lambda_{i_{\sigma(l)}})
\prod_{\substack{k<l\\\sigma(k)<\sigma(l)}}-h(\lambda_{i_{\sigma(k)}}-\lambda_{i_{\sigma(l)}})^{-1}\\
  &\qquad\qquad \times \prod_{\substack{k<l\\ \rho(k)<\rho(l)}}-h(\mu_{j_{\rho(k)}}- \mu_{j_{\rho(l)}})
t_{i_{\sigma(1)}j_{\rho(1)}}\cdots t_{i_{\sigma(r)}j_{\rho(r)}}\\
&\quad =\left(\sum_{\sigma\in S_r}
\prod_{\substack{k<l}}-h(\lambda_{i_{\sigma(k)}}-\lambda_{i_{\sigma(l)}})\right)
\eta^I_J.
\end{split}
\end{equation*}
So it suffices to prove that $A(I)(\ulambda):=\sum_{\rho\in
S_r}\prod_{k<l} -h(\lambda_{i_{\rho(k)}}-\lambda_{i_{\rho(l)}})$ is
independent of $\lambda$ and $I$:
\begin{equation*}
A(I)(\ulambda)=\sum_{\rho\in S_r}
\prod_{k<l}(-q)\frac{q^{-2\lambda_{i_{\rho(k)}}}-q^{-2}q^{-2\lambda_{i_{\rho(l)}}}}
{q^{-2\lambda_{i_{\rho(k)}}}-q^{-2\lambda_{i_{\rho(l)}}}}=
(-q)^{\frac{1}{2}r(r-1)}\prod_{k=1}^r \frac{1-q^{-2k}}{1-q^{-2}}\neq
0,
\end{equation*}
using the explicit expression \eqref{eq:defnh} for $h$ and
\eqref{eq:HLidentity}.
\end{proof}

\begin{cor} $\det=\widetilde{\det}$.
\end{cor}

\begin{remark}
\begin{mylist}
\item The dynamical quantum minor determinant $\xi^I_J$ belongs to
the weight space $\F_R(M(n))_{\omega(I),\omega(J)}$, where
$\omega(I)=\sum_{k=1}^r \omega(i_k)=\sum_{k=1}^r e_{i_k}$, $I=\{i_1,
\ldots, i_r\}$.
\item From Theorem \ref{thm:xi=eta} we obtain relations in
$\F_R(M(n))$. For $r=2$ we get quadratic relations for the
generators $t_{ij}$, for $\rho=\id$ in the expressions of $\xi^I_J$
and $\eta^I_J$ in \eqref{eq:xi} and \eqref{eq:eta} we get 
the third relation of \eqref{eq:RLLdynMn}. Similarly, from
Proposition \ref{prop:matrixelementsR} for $r=2$ and taking the
expressions for the dynamical quantum minor for $\rho=\id$ and
$\rho=(12)$ we get \eqref{eq:relationproofR}.
\end{mylist}
\end{remark}

\section{Laplace expansions}\label{sec:Laplace}

In this section, we prove some expansion formulas for the dynamical
quantum minor determinants, which are used in the following
section to introduce the antipode. 

For $I_1$, $I_2$ disjoint ordered subsets of $\{1,\ldots,n\}$,
denote by $\mathrm{sign}(I_1;I_2)$ the element of $\Mh$ defined by
\begin{equation*}
\mathrm{sign}(I_1;I_2)(\ulambda)=\prod_{\substack{k>m\\ k\in I_1,
m\in I_2}} -h(\lambda_k-\lambda_m).
\end{equation*}
Then $w_{I_1}w_{I_2}=\mu_W(\mathrm{sign}(I_1;I_2))w_I$ if $I_1\cap
I_2=\emptyset$ and $I_1\cup I_2=I$. If $I_1\cap I_2\neq\emptyset$
then $w_{I_1} w_{I_2}=0$ and in this case we define
$\mathrm{sign}(I_1;I_2)(\ulambda)=0$. For $I_1\cap I_2=\emptyset$
and $I=I_1\cup I_2$ as ordered subset we have
$\sign(I_1;I_2)=S(\sigma,I)$ where $\sigma$ is the permutation which
maps $I_1\cup I_2$ to the ordered subset $I$.

\begin{prop}[Laplace expansions]\label{prop:laplace}
Let $I$, $J_1$, $J_2$ be subsets of $\{1, \ldots , n\}$. If
$J=J_1\cup J_2$, $\# J=\# I$ we have
\begin{equation}\label{eq:laplace}
\begin{split}
\mu_r(\mathrm{sign}(J_1;J_2)) \xi^I_J &= \sum_{I_1 \cup I_2=I}
\mu_l(\mathrm{sign}(I_1;I_2))\xi_{J_1}^{I_1} \xi_{J_2}^{I_2},\\
\mu_l(T_{-\omega(J_1)}\mathrm{sign}(J_2;J_1)^{-1}) \xi^J_I &=
\sum_{I_1 \cup
I_2=I}\mu_r(T_{-\omega(I_1)}\mathrm{sign}(I_2;I_1)^{-1})
\xi^{J_1}_{I_1} \xi^{J_2}_{I_2},
\end{split}
\end{equation}
where the summation runs over all partitions $I_1 \cup I_2=I$ of $I$
such that $\#I_1=\#J_1$, $\#I_2=\#J_2$.
\end{prop}

\begin{remark}
\begin{mylist}
\item Note that the left hand sides of the expressions in
\eqref{eq:laplace} are zero if $J_1$ and $J_2$ are not disjoint.
\item The second relation of \eqref{eq:laplace} can be rewritten as
\begin{equation*}
\xi^J_I = \sum_{I_1 \cup I_2=I} \xi^{J_1}_{I_1}
\frac{\mu_l(\mathrm{sign}(J_2;J_1))}{\mu_r(\mathrm{sign}(I_2;I_1))}
\xi^{J_2}_{I_2}.
\end{equation*}
\end{mylist}
\end{remark}

\begin{proof}[Proof of Proposition \ref{prop:laplace}]
We have
\begin{equation*}
\begin{split}
R(w_{J_1})R(w_{J_2})
 & = \sum_{I_1\cap I_2= \emptyset}w_{I_1}w_{I_2} \otimes \xi^{I_1}_{J_1}\xi^{I_2}_{J_2}
   = \sum_{I_1\cap I_2= \emptyset} \mu_W(\mathrm{sign}(I_1;I_2)) w_I \otimes
\xi^{I_1}_{J_1}\xi^{I_2}_{J_2}\\
 & = \sum_{I_1\cap I_2= \emptyset} w_I
\otimes \mu_l(\mathrm{sign}(I_1;I_2))\xi_{J_1}^{I_1}
\xi_{J_2}^{I_2},
\end{split}
\end{equation*}
Also, if $J_1\cap J_2\neq \emptyset$ then $R(w_{J_1}w_{J_2})=R(0)=0$
by \eqref{eq:relationsW} which proves the first relation of
\eqref{eq:laplace} using Lemma \ref{lem:Wbasis} in the case that
$J_1$ and $J_2$ are not disjoint. If $J_1\cap J_2=\emptyset$ then we
also have
\begin{equation*}
\begin{split}
R(w_{J_1})R(w_{J_2})&= R(w_{J_1}w_{J_2})= (1\otimes
\mu_r(\mathrm{sign}(J_1;J_2)))R(w_J)  = \sum_{\#I=\#J} w_I \otimes
\mu_r(\mathrm{sign}(J_1;J_2))\xi^I_J.
\end{split}
\end{equation*}
The second relation of
\eqref{eq:laplace} is proved analogously, using $L$
instead of $R$ and Theorem \ref{thm:xi=eta}.
\end{proof}

In the special case $\#I=\#J=n$ and either $J_1$ or $J_2$ contains
one element, we get the following expansion formulas for the
determinant element. These expansions can be seen as dynamical
equivalent of the cofactor expansion across a row or column of the
determinant of a matrix.

\begin{cor}\label{cor:detontw}
For all $1\leq i,j \leq n$ we have
\begin{equation*}
\begin{split}
 \delta_{ij} \det = \sum_{k=1}^n \frac{\sign(\{k\};\hat k)(\ulambda)}{\sign(\{i\};\hat\imath)(\umu)}
      t_{kj}\xi_{\hat\imath}^{\hat k},& \qquad
  \delta_{ij} \det = \sum_{k=1}^n t_{jk} \frac{\sign(\hat\imath;\{i\})(\ulambda)}{\sign(\hat k;\{k\})(\umu)}
        \xi_{\hat k}^{\hat \imath}, \\
  \delta_{ij} \det = \sum_{k=1}^n \frac{\sign(\hat k;\{k\})(\ulambda)}{\sign(\hat\imath;\{i\})(\umu)}
        \xi_{\hat \imath}^{\hat k}t_{kj},& \qquad
  \delta_{ij} \det = \sum_{k=1}^n \xi_{\hat k}^{\hat\imath}
        \frac{\sign(\{i\};\hat\imath)(\ulambda)}{\sign(\{k\};\hat k)(\umu)}t_{jk},
\end{split}
\end{equation*}
with the notation $\hat \imath= \{1, \ldots, i-1, i+1, \ldots n\}$.
\end{cor}


\section{The dynamical $GL(n)$ quantum group}\label{sec:dynGLn}

In this section we extend $\F_R(M(n))$ by adjoining an
inverse of
the determinant. The resulting $\h$-bialgebroid $\F_R(GL(n))$ 
is equipped with an antipode, so it is 
a $\h$-Hopf algebroid. 

\begin{lemma}\label{lem:detcomm}
In $\F_R(M(n))$, the determinant element commutes with all quantum
minor determinants $\xi^{I}_{J}$, for $I$, $J$ subsets of
$\{1,\ldots n\}$. In particular, $\det$ commutes with all generators
$t_{ij}$.
 Moreover, $\Delta(\det)=\det\otimes \det$ and
$\varepsilon(\det)=T_{-\underline 1}$, with $\underline 1=
(1,\ldots,1)\in \h^*$.
\end{lemma}

\begin{proof}
Denote by $T$ the $n\times n$-matrix with elements $t_{ij}$, where
$i$ indicates the row index. Using the notation
\begin{equation}\label{eqn:notationT}
T_{j}^{i}= \frac{\mu_l(\sign(\hat \imath;\{i\}))}
{\mu_r(\sign(\hat\jmath;\{j\}))}\xi_{\hat \jmath}^{\hat \imath},
\end{equation}
denote by $\tilde T$ the $n\times n$-matrix with elements $T_j^i$
where $i$ indicates the column index. Then the third relation of
Corollary \ref{cor:detontw} implies $\tilde T T=\det\:I$ 
as $n^2$ identities in $\F_R(M(n))$, where $I$
is the $n\times n$-identity matrix. 
So $\det\: T= T\tilde T T=T\:\det $ which implies that $\det$
commutes with all generators $t_{ij}$. Since $\det\in
\F_R(M(n))_{\underline 1,\underline 1}$, we see that
$\det$ commutes with all elements in $M_{\h^\ast}$ that
only depend on differences $\lambda_i-\lambda_j$.  
By \eqref{eq:xi}, $\det$ also
commutes with $\xi^{I}_{J}$ for all subsets $I$, $J$. The last
statements follow from Corollary \ref{cor:comultxi}.
\end{proof}

So the determinant element commutes with all generators $t_{ij}$, but
since $\det\in \F_R(M(n))_{\underline 1,\underline 1}$ the
element $\det$ is not central. However, the set $S=\{ \det^k\}_{k\geq 1}$
satisfies the Ore condition, and this implies that we can 
localize at $\det$, see \cite{McCoR}. 
We  adjoin
$\F_R(M(n))$ with the formal inverse $\det^{-1}$, adding the
relations $\det\:\det^{-1}=1=\det^{-1}\det$,
$t_{ij}\det^{-1}=\det^{-1}t_{ij}$ and
$f(\ulambda)\det^{-1}=\det^{-1}f(\ulambda-\underline 1)$,
$f(\umu)\det^{-1}=\det^{-1}f(\umu-\underline 1)$. We denote the
resulting algebra by $\F_R(GL(n))$ and equip it with a bigrading
$\F_R(GL(n))=\bigoplus_{m,p\in\Z^n} \F_R(GL(n))_{mp}$ by
$\det^{-1}\in (\F_R(GL(n)))_{-\underline 1,-\underline 1}$. Lemma \ref{lem:detcomm} implies that $\det^{-1}$
commutes with all dynamical quantum minor determinants
$\xi^I_J$. By extending the comultiplication and counit of
Definition \ref{def:dynMn} by $\Delta (\det^{-1}) =
\det^{-1}\otimes\det^{-1}$, $\varepsilon (\det^{-1})=T_{\underline
1}$, $\F_R(GL(n))$ it is easily checked that 
$\F_R(GL(n))$ is a $\h$-bialgebroid.

\begin{prop}\label{prop:dynMnantipode}
The $\h$-bialgebroid $\F_R(GL(n))$ is a $\h$-Hopf algebroid with the
antipode $S$ defined on the generators by $S(\det^{-1})=\det$,
$S(\mu_r(f))= \mu_l(f)$, $S(\mu_l(f))=\mu_r(f)$ for all $f\in\Mh$
and
\begin{equation}\label{eq:dynMnantipode}
 \begin{split}
& S(t_{ij})=\det^{-1}
\frac{\mu_l(\sign(\hat\jmath;\{j\}))}{\mu_r(\sign(\hat\imath;\{i\}))}
\xi^{\hat\jmath}_{\hat\imath},
 \end{split}
\end{equation}
and extended as an algebra anti-homomorphism.
\end{prop}

\begin{proof} By \cite[Prop. 2.2]{KoelinkRosengren} it suffices to
check that $S$ is well-defined and
that \eqref{eq:defS} holds on the generators. It is straightforward
to check that $S$ preserves the relations \eqref{eq:commdynMnftij}. 
To see that $S$ preserves the RLL-relations, we apply
the antipode to the RLL-relations \eqref{eq:RLL}. 
Using \eqref{eqn:notationT} this gives 
\begin{equation}\label{eq:SRLLdirect}
\sum_{x,y} \det^{-2}T^d_y T^b_x R^{xy}_{ac}(\umu)=\sum_{x,y}
\det^{-2} T^x_a T^y_c R^{bd}_{xy}(\ulambda),
\end{equation}
which is equivalent to
\begin{equation}\label{eq:SRTT}
\sum_{x,y} R_{ac}^{xy}(\umu+\omega(x)+\omega(y))T_y^d T_x^b
\det^{-2}= \sum_{x,y} \det^{-2}T_a^x T_c^y R_{xy}^{bd}(\ulambda).
\end{equation}
We have to prove that \eqref{eq:SRTT} holds in $\F_R(GL(n))$. To
show this, we multiply the RLL-relations \eqref{eq:RLL} by $T^k_d
T_b^l$ from the right and by $ T^a_j T^c_i$ from the left and sum
over all $a$, $b$, $c$ and $d$ we get, using Corollary
\ref{cor:detontw},
\begin{equation}\label{eq:commTT}
\sum_{a,c}T^a_j T^c_i R^{lk}_{ac}(\ulambda)\det^2= \sum_{b,d}\det^2
R^{bd}_{ji}(\umu+\omega(i)+\omega(j))T^k_d T^l_b,
\end{equation}
Multiplying this equation from the left and from the right by
$\det^{-2}$ gives \eqref{eq:SRTT} by the $\h$-invariance
of the $R$-matrix, so $S$ preserves the
RLL-relations. 

{}From the proof of Lemma \ref{lem:detcomm} it follows that $S(T)
T=TS(T)=I$, where $T$ is defined as in the proof of Lemma
\ref{lem:detcomm}, so \eqref{eq:defS} holds for all generators
$t_{ij}$. The proof of \cite[Prop. 2.2]{KoelinkRosengren} shows that  if \eqref{eq:defS} holds for $a$ and $b$, then it holds for $ab$, so that
in particular \eqref{eq:defS} holds for $\det$. By Lemma
\ref{lem:detcomm} we find $S(\det)\det=1=\det S(\det)$, so that
$S(\det)=\det^{-1}$. An independent proof of this statement
is given in Proposition \ref{prop:Sonquantuminordet}. 
With this observation it is easily proved that $S$
also preserves the defining relations involving $\det^{-1}$, and
that \eqref{eq:defS} holds for $\det^{-1}$. 
\end{proof}

The relation $S(\det)=\det^{-1}$ is the special case 
$I=J=\{1,\ldots,n\}$ of 
the following proposition. 

\begin{prop}\label{prop:Sonquantuminordet} 
For $I$ and $J$ ordered subsets such that
$\# I=\# J$ we have
\begin{equation}\label{eq:dynMnSonxi}
S(\xi^I_J)=\det^{-1}\frac{\mu_l(\sign(J^c;J))}{\mu_r(\sign(I^c;I))}
\xi^{J^c}_{I^c},
\end{equation}
with $I^c$ the complement of $I$ in $\{1,\ldots,n\}$.
\end{prop}

\begin{proof}
We prove this formula by induction on the $r:=\# I$ using the Laplace expansions of Proposition \ref{prop:laplace}. Another proof uses
\eqref{eq:defS} combined with the Laplace expansions. 
We use a similar induction step in the proof of 
Lemma \ref{lem:dynGLn*structure}. 

For $r=1$ this is just the definition of the antipode on a
generator $t_{ij}$. For the induction step we use the Laplace
expansions of Proposition \ref{prop:laplace} several times. Let
$j\in J$, then applying the Laplace expansion twice
\begin{equation*}
\begin{split}
S(\xi^I_J)& = S\left( \sum_{i\in I}
 \frac{\mu_l(\sign(\{i\};I'))}{\mu_r(\sign(\{j\};J'))}\xi^{\{i\}}_{\{j\}}
\xi^{I'}_{J'}\right) =\sum_{i\in I} S(\xi^{I'}_{J'})S(\xi^i_j)
\frac{\mu_r(\sign(\{i\};I'))}{\mu_l(\sign(\{j\};J'))}\\
& = \sum_{i\in I} \det^{-2}
\frac{\mu_l(\sign(J'^c;J'))}{\mu_r(\sign(I'^c;I'))}\xi^{J'^c}_{I'^c}
 \frac{\mu_l(\sign(\hat\jmath;\{j\}))}{\mu_r(\sign(\hat\imath;\{i\}))}\xi^{\hat\jmath}_{\hat\imath}
 \frac{\mu_r(\sign(\{i\};I'))}{\mu_l(\sign(\{j\};J'))}\\
&= \sum_{i\in I} \det^{-2}
\frac{\mu_l(\sign(J^c;J'))}{\mu_r(\sign(I^c;I'))}\xi^{J'^c}_{I'^c}
 \frac{\mu_l(\sign(\hat\jmath;\{j\}))}{\mu_r(\sign(\hat\imath;\{i\}))}\xi^{\hat\jmath}_{\hat\imath}\\
&=\sum_{i\in I} \sum_{\substack{k\in J'^c\\K=J'^c\setminus\{k\}}}
 \det^{-2}\frac{\mu_l(\sign(J^c;J'))}{\mu_r(\sign(I^c;I'))}
 \frac{\mu_l(\sign(K;\{k\}))}{\mu_r(\sign(I^c;\{i\}))}\xi^K_{I^c}
 \xi^{\{k\}}_{\{i\}}
 \frac{\mu_l(\sign(\hat\jmath;\{j\}))}{\mu_r(\sign(\hat\imath;\{i\}))}\xi^{\hat\jmath}_{\hat\imath},
\end{split}
\end{equation*}
where $J'= J\setminus \{j\}$ and $I'=I\setminus\{i\}$ (so $I'$
depends on the summation index) as ordered subsets. In this
computation we use $\xi^{J'^c}_{I'^c}\xi^{\hat\jmath}_{\hat\imath}
\frac{\mu_r(\sign(\{i\};I'))}{\mu_l(\sign(\{j\};J'))}=\frac{\mu_r(\sign(\{i\};I'))}{\mu_l(\sign(\{j\};J'))}
\xi^{J'^c}_{I'^c}\xi^{\hat\jmath}_{\hat\imath}$ since
$\xi^{J'^c}_{I'^c}\xi^{\hat\jmath}_{\hat\imath}$ has weight
$(\underline 1+\omega(J^c), \underline 1 +\omega(I^c))$ and
$\sign(\{a\};B)\sign(A;B)=\sign(\{a\}\cup A;B)$ for all subsets $A$,
$B$ and all elements $a\not\in A$. Since $\sum_{\substack{k\in
J'^c\\K=J'^c\setminus\{k\}}}
\mu_l(\sign(K;\{k\}))\xi^K_{I^c}\xi^{\{k\}}_{\{i\}}=0$ for all $i\in
I^c$ and $\sign(\hat\imath;\{i\})=\sign(I';\{i\})\sign(I^c;\{i\})$
we obtain, using the Laplace expansion once more for the summation
over $i$ where the only non-zero term is for $k=j$,
\begin{equation*}
\begin{split}
S(\xi^I_J) &= \sum_{\substack{k\in J'^c\\K=J'^c\setminus\{k\}}}
\det^{-2}\frac{\mu_l(\sign(J^c;J'))}{\mu_r(\sign(I^c;I))}
 {\mu_l(\sign(K;\{k\}))}\xi^K_{I^c} \sum_{i=1}^n \xi^{\{k\}}_{\{i\}} \frac{\mu_l(\sign(\hat\jmath;\{j\}))}
{\mu_r(\sign(\hat\imath;\{i\}))}
\xi^{\hat\jmath}_{\hat\imath}\\
 &=\det^{-2}\frac{\mu_l(\sign(J^c;J'))}{\mu_r(\sign(I^c;I))}
 {\mu_l(\sign(J^c;\{j\}))}\xi^{J^c}_{I^c}\det =\det^{-1}\frac{\mu_l(\sign(J^c;J))}{\mu_r(\sign(I^c;I))}\xi^{J^c}_{I^c},
\end{split}
\end{equation*}
which proves the proposition.
\end{proof}

\begin{cor}\label{cor:Ssquaredonquantumminor}
\begin{equation*}
\begin{split}
S^2(\xi^I_J)& =
 \frac{\prod_{m\in I,k\in I^c} h(\lambda_m-\lambda_k)}{\prod_{m\in
J,k\in J^c} h(\mu_m-\mu_k)}\xi^I_J.
\end{split}
\end{equation*}
\end{cor}

In particular, $S$ is invertible.

\section{The dynamical $U(n)$ quantum group}\label{sec:dynUn}

In this section we prove the existence of a $\ast$-operator 
on $\F_R(GL(n))$, such that it becomes a $\h$-Hopf $\ast$-algebroid.
Equipped with this $\ast$-structure we denote the 
$\h$-Hopf $\ast$-algebroid by $\F_R(U(n))$. 

\begin{lemma}\label{lem:dynGLn*structure}
The $*$-operator defined on
the generators by
\begin{equation*}
t_{ij}^*=\xi^{\hat \imath}_{\hat \jmath} \det^{-1},\qquad
\mu_l(f)^*=\mu_l(\overline f),\; \mu_r(f)=\mu_r(\overline f), \qquad
(\det^{-1})^*=\det,
\end{equation*}
and extended as $\C$-antilinear algebra anti-homomorphism 
is well-defined on $\F_R(GL(n))$. 
\end{lemma}

\begin{proof}
Let $I$ and $J$ be ordered subsets of $\{1, \ldots, n\}$, such that
$\#I= \#J=r$. Denote by $I^c$ the complement of $I$ in $\{1, \ldots,
n\}$, then we have
\begin{equation}\label{eq:dynMn*onxi}
(\xi^I_J)^*= \xi^{I^c}_{J^c} \det^{-1}.
\end{equation}
{}From this result and Lemma \ref{lem:detcomm} it directly follows
that $*$ is an involution. The proof of \eqref{eq:dynMn*onxi} is
analogous to the corresponding statement \eqref{eq:dynMnSonxi} for
the antipode.

We prove that $*$ preserves the RLL-relations by using that the
antipode does so. By definition of $S$ and $*$ it follows that
\begin{equation}
\mu_r(\sign(\hat k;\{k\}))S(t_{kj})=
\mu_l(\sign(\hat\jmath;\{j\}))t_{jk}^*.
\end{equation}
Applying $*$ to the RLL-relations \eqref{eq:RLL} we get
\begin{equation*}
\begin{split}
\sum_{x,y=1}^n\frac{\mu_r(\sign(\hat d;\{d\}))}{\mu_l(\sign(\hat
y;\{y\}))}&T^y_d \frac{\mu_r(\sign(\hat b;\{b\}))}{\mu_l(\sign(\hat
x;\{x\}))}T^x_b \mu_l(R^{xy}_{ac})\\
 &= \sum_{x,y=1}^n \frac{\mu_r(\sign(\hat x;\{x\}))}{\mu_l(\sign(\hat
a;\{a\}))}T^a_x \frac{\mu_r(\sign(\hat y;\{y\}))}{\mu_l(\sign(\hat
c;\{c\}))}T^c_y \mu_r(R^{bd}_{xy}),
\end{split}
\end{equation*}
which is equivalent to
\begin{equation*}
\begin{split}
\sum_{x,y=1}^n\frac{\mu_l(\sign(\hat a;\{a\}))}{\mu_l(\sign(\hat
y;\{y\}))} &\frac{\mu_l(T_{\omega(a)}\sign(\hat
c;\{c\}))}{\mu_l(T_{\omega(y)}\sign(\hat x;\{x\}))}T^y_d T^x_b
\mu_l(R^{xy}_{ac})\\
 &= \sum_{x,y=1}^n \frac{\mu_r(\sign(\hat x;\{x\}))}{\mu_r(\sign(\hat d;\{d\}))}
\frac{\mu_r(T_{\omega(x)}\sign(\hat y;\{y\}))}
{\mu_r(T_{\omega(d)}\sign(\hat b;\{b\}))} T^a_x T^c_y
\mu_r(R^{bd}_{xy}).
\end{split}
\end{equation*}
Using \eqref{eq:SRLLdirect}, $*$ preserves the RLL-relations if
\begin{equation}\label{eqn:relationelementsR}
R^{yx}_{db}(\umu)=R^{bd}_{xy}(\umu)\frac{\sign(\hat
x;\{x\})(\umu-\omega(x)-\omega(y))}{\sign(\hat
d;\{d\})(\umu-\omega(x)-\omega(y))} \frac{\sign(\hat
y;\{y\})(\umu-\omega(y))} {\sign(\hat b;\{b\})(\umu-\omega(b))}.
\end{equation}
This follows by direct calculations using the explicit expression of
$R$ and the fact that $\sign(\hat x;x)$ is independent of $\mu_y$
for all $y<x$, where the only non-trivial cases are for $x=y=b=d$,
$x=b,y=d$ and $x=d,y=b$. Using $\det^*=\det^{-1}$ which follows from
\eqref{eq:dynMn*onxi}, it directly follows that $*$ preserves the
other commutation relations.
\end{proof}

\begin{prop}\label{prop:FRUn}
Denote $\F_R(GL(n))$ equipped with the $\ast$-operator
of Lemma \ref{lem:dynGLn*structure} by $\F_R(U(n))$,
then $\F_R(U(n))$ is a $\h$-Hopf $\ast$-algebroid.
\end{prop}

\begin{proof}
From the definition of $\ast$ and Corollary \ref{cor:comultxi} it follows that
$(\ast\otimes \ast)\Delta (t_{ij})=\Delta(t_{ij}^\ast)$ and
$(\varepsilon\circ \ast)(t_{ij}) =(\ast^{\Dh}\circ \varepsilon)(t_{ij})$,
\begin{equation*}
\begin{split}
 &(*\otimes *)\Delta(\det^{-1})=\det\otimes\det=\Delta((\det^{-1})^*),\qquad
  (\varepsilon\circ *)(\det^{-1})= T_{-\underline 1}= (T_{\underline
1})^*=\varepsilon(\det^{-1})^*.
\end{split}
\end{equation*}
So the relations $(*\otimes *)\circ \Delta=\Delta\circ *$ and
$\varepsilon\circ *=*^{\Dh}\circ \varepsilon$ hold on the generators
of $\F_R(GL(n))$ and hence on all of $\F_R(GL(n))$. 
\end{proof}

From \eqref{eq:dynMnSonxi} and \eqref{eq:dynMn*onxi} it directly
follows that
\begin{equation}\label{eq:*SandS*}
\begin{split}
&
S(\xi^I_J)^*=\xi^J_I\frac{\mu_l(\sign(J^c;J))}{\mu_r(\sign(I^c;I))},\qquad
S((\xi^I_J)^*)=\frac{\mu_l(\sign(J;J^c))}{\mu_r(\sign(I;I^c))}\xi^J_I,
\end{split}
\end{equation}
which gives an indication for the unitarisability of the
corepresentations $R$ and $L$ of $\F_R(GL(n))$ defined in
Proposition \ref{prop:extensionR} and \ref{prop:extensionL}, for the
definition of unitarisability see \cite[\S5]{KvN}.

\begin{prop}
The corepresentations $R$ and $L$ are unitarisable corepresentations
of $\F_R(U(n))$.
\end{prop}

\begin{proof}
We have to define a form $\pair{\cdot}{\cdot}:W\times W\to \Mh$ and
check that $\pair{R(x)}{R(y)}=\mu_r(\pair{x}{y}\1)$ for all $x$,
$y\in W$, see \cite[\S 5]{KvN}. 
It is sufficient to do this for basis elements $\{w_I\}$
of $W$. Define
$\pair{w_I}{w_J}(\ulambda)=\delta_{IJ}\sign(I^c;I)(\ulambda-\omega(I))\in
\Mh$, so $\pair{w_I}{w_J}_D=\delta_{IJ}\sign(I^c;I)\in\Dh$. Then
\begin{equation*}
\begin{split}
\pair{R(w_I)}{R(w_J)}&=\pair{\sum_{\# K=\# J} w_K\otimes
\xi^K_J}{\sum_{\#M =\#J} w_M\otimes \xi^M_J}=\sum_{K,M}
\pair{w_K}{w_M}_D\otimes
(\xi^M_J)^*\xi^K_J\\
& =\sum_K \mu_l(\sign(K^c;K))\frac{\mu_r(\sign(J^c;J))}
{\mu_l(\sign(K^c;K))} S(\xi^J_K)\xi^K_I\\
 &= \mu_r(\sign(J^c;J))\delta_{IJ} = \mu_r(\pair{w_I}{w_J}_D\1),
\end{split}
\end{equation*}
using \eqref{eq:*SandS*} and \eqref{eq:defS} on $\xi^J_I$. Define a
form on $V$ by $\pair{v_I}{v_J}=\delta_{IJ}\sign(I;I^c)^{-1}\in
\Mh$. By a similar computation it follows that
$\pair{L(v_I)}{L(v_J)}=\mu_l(\pair{v_I}{v_J}_D\1)$.
\end{proof}

\begin{remark}\label{rem:subalgebras}
The above discussion strongly suggests that there
are analogues of the dynamical $SL(n)$ and $SU(n)$
quantum groups. We refer to \cite{vNord} for details. 
\end{remark}


\section{A pairing on the dynamical $U(n)$ quantum group}
\label{sec:pairing}

In this section we discuss pairings for the dynamical
$GL(n)$ quantum group and we  present a cobraiding on $\F_R(GL(n))$. For
a pairing for $\F_R(GL(n))^\cop$ and $\F_R(GL(n))$  as $\h$-Hopf
$\ast$-algebroids, we need a second $*$-operator on $\F_R(GL(n))$.

\subsection{Pairing for $\h$-Hopf $*$-algebroids}

We start by recalling the definition of a pairing  for $\h$-Hopf $\ast$-algebroids.

\begin{defn}\label{def:pairingforhbialgebras}
A pairing for $\h$-bialgebroids $\U$ and $\A$ is a $\C$-bilinear map
$\langle \cdot,\cdot\rangle \colon \U\times \A \to \Dh$ satisfying
\begin{subequations}\label{eq:pairingforhbialgebras}
\begin{align}
\label{eq:pairing-target} &
\pair{\U_{\alpha\beta}}{\A_{\gamma\delta}}
\subseteq (\Dh)_{\alpha+\delta,\beta+\gamma}, \\
\label{eq:pairing-leftmoment} & \pair{\mu_l^\U(f)X}{a} =
\pair{X}{\mu_l^\A(f) a}
 = f \circ \pair{X}{a},\qquad
 \pair{X\mu_r^\U(f)}{a}= \pair{X}{a\mu_r^\A(f)}
 = \pair{X}{a} \circ f, \\
\label{eq:pairing-leftproduct} &\pair{XY}{a}=\sum_{(a)}
\pair{X}{a_{(1)}}T_\rho \pair{Y}{a_{(2)}},
    \quad \Delta^\A(a)=\sum_{(a)} a_{(1)}\otimes a_{(2)},\;
    a_{(1)}\in \A_{\gamma\rho},\\
\label{eq:pairing-rightproduct} &\pair{X}{ab}=\sum_{(X)}
\pair{X_{(1)}}{a}T_\rho\pair{X_{(2)}}{b},
     \quad \Delta^\U(X)=\sum_{(X)} X_{(1)}\otimes X_{(2)},\;
     X_{(1)}\in \U_{\alpha\rho},\\
\label{eq:pairing-counit-unit} &\pair{X}{1}=\varepsilon^\U(X),
\qquad \pair{1}{a}=\varepsilon^\A(a),
\end{align}
\end{subequations}
for all $X\in\U$, $a\in\A$. If moreover, $\U$ and $\A$ are $\h$-Hopf
algebroids, then in addition we require
\begin{equation}\label{eq:pairingforhHopfalgebras}
\pair{S^\U(X)}{a} = S^{\Dh}( \pair{X}{S^\A(a)}), \mbox{ for all
$X\in\U$, $a\in\A$.}
\end{equation}
If in addition a $*$-operator is defined on $\U$ and $\A$ such that
\begin{equation}\label{eq:pairingXstara}
\pair{X^\ast}{a} = T_{-\gamma}\circ
(\pair{X}{S^\A(a)^\ast})^\ast\circ T_{-\delta}, \mbox{ for all
$a\in\A_{\gamma\delta}$, $X\in\U$},
\end{equation}
then $\U$ and $\A$ are paired as $\h$-Hopf $*$-algebroids.
\end{defn}

\begin{remark}\label{rem:pairingoftenzero}
Note that \eqref{eq:pairing-target}
implies that $\pair{X}{a}=0$ whenever $X\in\U_{\alpha\beta}$,
$a\in\A_{\gamma\delta}$ with $\alpha+\delta\not=\beta+\gamma$.
\end{remark}

A cobraiding on a $\h$-bialgebroid $\A$ is a pairing
$\pair{\cdot}{\cdot}:\A^\cop \times \A\to \Dh$ which in addition
satisfies
\begin{equation}\label{eq:defcobraiding}
\sum_{(a),(b)} \mu_l^\A(\pair{a_{(1)}}{b_{(1)}}\1) a_{(2)}b_{(2)} =
\sum_{(a),(b)} \mu_r^\A(\pair{a_{(2)}}{b_{(2)}}\1) b_{(1)}a_{(1)},
\end{equation}
as an identity in $\A$ and where $\Delta^\A(a)=\sum_{(a)}
a_{(1)}\otimes a_{(2)}$, $\Delta^\A(b)=\sum_{(b)} b_{(1)}\otimes
b_{(2)}$. In \cite{r}, Rosengren proved that for a $\h$-bialgebroid
constructed by the generalized FRST-construction from an $R$-matrix,
denoted by $R$, that satisfies the quantum dynamical Yang-Baxter
equation \eqref{eq:QDYBE} there exists a natural cobraiding defined
on the generators by
\begin{equation}\label{eqn:cobrFRST}
\pair{L_{ij}}{L_{kl}}=R_{ik}^{jl}(\lambda)T_{-\omega(i)-\omega(k)}.
\end{equation}
Note that this is the dynamical analogue of the cobraiding
for quantum groups, see e.g. \cite[\S VIII.6]{Kass}

In \cite{KvN} we proved the following proposition, which we now
extend to the level of $\h$-(co)module algebras. By $\A^\lr$ we
denote the $\h$-algebra obtained from a $\h$-algebra $\A$ by
interchanging the moment maps and with weight spaces
$(\A^\lr)_{\alpha\beta}=\A_{\beta\alpha}$.

\begin{prop}\label{prop:cotodynrep}
Let $\U$ be a $\h$-algebra and $\A$ be $\h$-coalgebroid  equipped
with a pairing $\pair{\cdot}{\cdot}:\U\times \A\to \Dh$, and let $V$
be a $\h$-space.
\begin{mylist}
\item \label{item:cotodynrepright} Let $R\colon V\to V\totimes \A$ be a right corepresentation of the
$\h$-coalgebroid $\A$, then $\pi(X)v = (\id\otimes
\pair{X}{\cdot}T_\beta)R(v)$ for $X\in\U_{\alpha\beta}$, defines a
$\h$-algebra homomorphism $\pi\colon \U \to (D_{\hs,V})^\lr$, hence
$\pi\colon \U^\lr \to D_{\hs,V}$ defines a dynamical representation
of $\U^\lr$ on $V$.
\item \label{item:cotodynrepleft} Let $L\colon V\to \A\totimes V$ be a left corepresentation of the
$\h$-coalgebroid $\A$, then $\pi(X)v =
(T_{\alpha}\pair{X}{\cdot}\otimes\id)L(v)$ for
$X\in\U_{\alpha\beta}$, defines a $\h$-algebra homomorphism
$\pi\colon \U^\opp \to (D_{\hs,V})^\lr$. In particular, $\pi\colon
(\U^\opp)^\lr \to D_{\hs,V}$ defines a dynamical representation of
$(\U^\opp)^\lr$ on $V$. Moreover, if $\U$ is $\h$-Hopf algebroid ,
then $X\mapsto \pi(S^\U(X))$ defines a dynamical representation of
$\U$ on $V$.
\end{mylist}
\end{prop}

We now extend this result to the level of $\h$-comodule algebras.

\begin{defn}\label{def:hmodule}
Let $\A$ be a $\h$-bialgebroid and $V$ a $\h$-space. We call $V$ a
$\h$-module algebra for $\A$ if there exists a dynamical
representation $\pi:\A\to D_{\h^*,V}$ such that
\begin{mylist}
\item \label{item:hmoduleV} $V$ is an associative algebra such that 
$\mu_V(f) vw = v \mu_V(T_\alpha f)w$
for all $v\in V_\alpha$, $w\in V$, and $V_\alpha V_\beta\subset V_{\alpha+\beta}$,
\item $\pi(a)vw =\sum_{(a)} (\pi(a_{(1)})v)(\pi(a_{(2)})w)$, for
all $v$, $w\in V$ and $X\in\A$ with $\Delta(a)=\sum_{(a)}
a_{(1)}\otimes a_{(2)}$.
\end{mylist}
Moreover, if $V$ is unital then $\pi(a)1=\mu_V(\varepsilon(a)\1)$.
\end{defn}

\begin{prop}\label{prop:hmodule}
Let $\U$ and $\A$ be paired as $\h$-bialgebroids. Let $V$ be a right
(left) $\h$-comodule algebra for $\A$, then $\pi$ as defined in
Proposition \ref{prop:cotodynrep} defines a $\h$-module algebra for
$\U^{\lr}$ ($(\U^\opp)^\lr$).
\end{prop}

\begin{proof}
We prove the proposition in the case that $V$ is a right
$\h$-comodule algebra, the other statement can be proved
analogously. Since $V$ is a $\h$-comodule algebra Definition
\ref{def:hmodule} \ref{item:hmoduleV} is satisfied. By Proposition
\ref{prop:cotodynrep}, $\pi(X)v=(\id\otimes \pair{X}{\cdot
}T_\beta)R(v)$, $X\in\U_{\alpha\beta}$, is a $\h$-algebra
homomorphism of $\U$ to $(D_{\h^*,V})^\lr$. Then, since $R$ is an
algebra homomorphism we have
\begin{equation*}
\begin{split}
\pi(X)vw&= (\id\otimes \pair{X}{ \cdot}T_\beta)R(vw)=\sum v_{(1)}w_{(1)}\otimes \pair{X}{a_{(2)}b_{(2)}}T_{\beta}\\
  &=\sum v_{(1)}w_{(1)}\otimes \pair{X_{(1)}}{a_{(2)}}T_{\gamma}
 \pair{X_{(2)}}{b_{(2)}}T_{\beta}=\sum
 (\pi(X_{(1)})v)(\pi(X_{(2)})w),
\end{split}
\end{equation*}
for $X\in\U_{\alpha\beta}$, $\Delta(X)=\sum_{(X)} X_{(1)}\otimes
X_{(2)}$, $X_{(1)}\in \U_{\alpha\gamma}$ and with the notation
$R(v)=\sum v_{(1)}\otimes a_{(2)}$, $R(w)=\sum w_{(1)}\otimes
b_{(2)}$. So $\pi$ defines a $\h$-module algebra for $\U^\lr$.

\par\noindent If $V$ is unital then $\pi(X)1=1\otimes \pair{X}{1}T_{\beta} =\mu_V(\varepsilon(X)\1)$
for $X\in\U_{\alpha\beta}$.
\end{proof}

\subsection{A pairing on the dynamical $GL(n)$ quantum group}

A natural cobraiding on the algebra $\F_R(M(n))$ is given by
\eqref{eqn:cobrFRST}. For this pairing we have
$\pair{t_{ij}}{\det}=\delta_{ij} q T_{-\underline 1 -\omega(i)}$.
For normalisation purposes we multiply the pairing of two generators 
with a factor $q^{-1/n}$.
So we use the pairing $\pair{\cdot}{\cdot}: \F_R(M(n))^\cop\times
\F_R(M(n))\to \Dh$ defined on the generators $t_{ij}$ by
\begin{equation}\label{eq:cobrdynMn}
\pair{t_{ij}}{t_{kl}}=q^{-1/n}R_{ik}^{jl}(\ulambda)T_{-\omega(i)-\omega(k)}.
\end{equation}
Note that switching $R$ to $q^{-1/n}R$ is a gauge transform, which
does not affect the RLL-relations. 
The non-trivial cases for this pairing on the level of the
generators are explicitly given by
\begin{equation}\label{eq:explicitpairingtij}
\begin{split}
 &\pair{t_{ii}}{t_{ii}} = q^{1-1/n} T_{-2\omega(i)}, \mbox{for all } i,\\
 &\pair{t_{ii}}{t_{jj}} = q^{-1/n}T_{-\omega(i)-\omega(j)}, \mbox{for all }i<j,\\
 &\pair{t_{jj}}{t_{ii}} = q^{-1/n}g(\lambda_i-\lambda_j) T_{-\omega(i)-\omega(j)}, \mbox{for all }i<j,\\
 &\pair{t_{ji}}{t_{ij}} = q^{-1/n}h_0(\lambda_i-\lambda_j)T_{-\omega(i)-\omega(j)}, \mbox{for all }i\neq j .
\end{split}
\end{equation}
In this section we prove that this pairing can be extended the level
of $\h$-Hopf $*$-algebroids.

\vspace{\baselineskip}

In order to extend the pairing to a cobraiding on $\F_R(GL(n))$ we
need to compute the pairing of a generator $t_{ij}$ with the
determinant element. Denote by $\underline 1$ the vector with all
$1$'s.

\begin{lemma}\label{lem:pairingtijdet}
For the pairing $\pair{\cdot}{\cdot}:\F_R(M(n))^\cop\times
\F_R(M(n))\to \Dh$ defined in \eqref{eq:cobrdynMn} we have
\begin{equation*}
 \pair{t_{ij}}{\det}= \delta_{ij} T_{-\underline 1-\omega(i)},\qquad
    \pair{\det}{t_{ij}}= \delta_{ij} T_{-\underline 1-\omega(i)},
    \qquad \pair{\det}{\det}=T_{-2\cdot \underline 1}.
\end{equation*}
\end{lemma}

\begin{proof}
{}From Remark \ref{rem:pairingoftenzero} it immediately follows that
$\pair{\det}{t_{ij}}=\pair{t_{ij}}{\det}=0$ for $i\neq j$. Using the
pairing \eqref{eq:cobrdynMn} on $\F_R(M(n))$, Propositions
\ref{prop:extensionR} and \ref{prop:hmodule} show that
$\pi:(\F_R(M_n)^\cop)^\lr \to (D_{\h^*,W})$ gives $W$ a $\h$-module
algebra structure for $(\F_R(M(n))^\cop)^\lr$. Then we have
\begin{equation*}
\pi(t_{ii})(w_1\cdots w_n)= w_1\cdots w_n \otimes
\pair{t_{ii}}{\det}T_{\omega(i)}.
\end{equation*}
Also we compute
\begin{equation*}
\begin{split}
\pi(t_{ii})w_1\cdots w_n  &=
\pi(t_{ii})\left(\prod_{k<i}-h(\lambda_i-\lambda_k)^{-1} w_i
w_{\hat\imath}\right)
 =
 \left(T_{-\omega(i)}\prod_{k<i}-h(\lambda_i-\lambda_k)^{-1}\right)
\pi(t_{ii}) [w_i w_{\hat\imath}]\\
 & =\prod_{k<i}-h(\lambda_i-1-\lambda_k)^{-1}\sum_{k_1,\ldots,k_n}
\sum_{j_1,\ldots,j_{n-1}} w_{k_1}w_{k_2}\cdots w_{k_n}\\
 &\qquad   \otimes
\pair{t_{j_1 i}}{t_{k_1 i}}T_{\omega(j_1)} \pair{t_{j_2 j_1}}{t_{k_2
1}} \cdots T_{\omega(j_{n-1})} \pair{t_{i j_{n-1}}}{t_{k_n
n}}T_{\omega(i)},
\end{split}
\end{equation*}
using the $\h$-module algebra structure of $W$ in the third
equation. From \eqref{eq:explicitpairingtij} it follows that
$\pair{t_{j_1i}}{t_{k_1i}}\neq 0$ only if $j_1=k_1=i$. So we get
\begin{equation*}
\begin{split}
\prod_{k<i}-h(\lambda_i-1-\lambda_k)^{-1}&\sum_{k_2,\ldots,k_n}
\sum_{j_2,\ldots,j_{n-1}} w_{i}w_{k_2}\cdots w_{k_n} \\ & \otimes
\pair{t_{i i}}{t_{i i}}T_{\omega(j_1)} \pair{t_{j_2 j_1}}{t_{k_2 1}}
\cdots T_{\omega(j_{n-1})} \pair{t_{i j_{n-1}}}{t_{k_n
n}}T_{\omega(i)}.
\end{split}
\end{equation*}

 Now $\pair{t_{j_2i}}{t_{k_21}}\neq 0$ only if $k_2=i$, $j_2=1$ or
$k_2=1$, $j_2=i$. In the first case, the first leg of the tensor
product is equal to $0$, so $k_2=1$, $j_2=i$. Continuing in this way
and recalling that we have pulled the term corresponding to $w_i$ to
the left, we obtain that there is only a non-zero contribution for
$j_m=i$ for all $m$ and $k_1=i$, $k_m=m-1$ for $2\leq m\leq i$ and
$k_m=m$ for $m>i$. So we get
\begin{equation*}
\begin{split}
 \pi(t_{ii})w_1&\cdots w_n  \\
  & =\prod_{k<i}-h(\lambda_i-1-\lambda_k)^{-1}
w_{i}w_{\hat\imath} \otimes
\pair{t_{ii}}{t_{ii}}T_{\omega(i)}\pair{t_{ii}}{t_{11}}T_{\omega(i)}
\cdots \pair{t_{ii}}{t_{nn}}T_{\omega(i)}\\
 & = \prod_{k<i}-h(\lambda_i-1-\lambda_k)^{-1}
\prod_{k<i}-h(\lambda_i-\lambda_k)w_1\cdots w_n \\
 & \quad \otimes q^{1-1/n}T_{-\omega(i)} \prod_{k<i}q^{-1/n}
g(\lambda_k-\lambda_i)T_{-\omega(k)}
\prod_{k>i}q^{-1/n}T_{-\omega(k)}\\
 &= \prod_{k<i}\frac{h(\lambda_i-\lambda_k)}{h(\lambda_i-\lambda_k-1)}
\prod_{k<i} g(\lambda_i-\lambda_k-1) w_1\cdots w_n =w_1\cdots w_n,
\end{split}
\end{equation*}
where the last equality follows from \eqref{eq:relationshg}. So
$\pair{t_{ii}}{\det}=T_{-\underline 1-\omega(i)}$.

Note that $\F_R(M(n))^\cop$ can also be seen as a $\h$-bialgebroid
constructed from the $R$-matrix $\tilde R$ with matrix elements
$\tilde R_{ab}^{cd}=R_{dc}^{ba}$ by the generalized
FRST-construction. Following the lines of the proofs of \S $3$ we
can prove that $V$ is a right $\h$-comodule algebra for
$\F_R(M(n))^\cop$. By inspection it follows that 
the matrix elements $\tau^I_J$
of this corepresentation $R^\cop$, defined by $R^\cop(v_I)=\sum_J
v_J\otimes \tau^I_J$, are equal to $\xi^I_J$. From Proposition
\ref{prop:hmodule} it follows that $\pi:\F_R(M(n))\to D_{\h^*,V}$
defined by $\pi(a)v=(\id \otimes \pair{\cdot}{a}T_{\beta})R^\cop(v)$
for $a\in \F_R(M_n)_{\alpha\beta} $ and $v \in V$ gives $V$ the
structure of a $\h$-module algebra for $\F_R(M(n))$. Now analogously
to the proof of the first part of this lemma we get
$\pair{\det}{t_{ii}}=T_{-\underline 1-\omega(i)}$.

Using Lemma \ref{lem:pairingtijdet}, $\pi(t_{ij})w_1\cdots w_n=0$ if
$i\neq j$ and the explicit expression of $\det$ we get
\begin{equation*}
\pi(\det)w_1\cdots w_n = \pi(t_{11}t_{22}\cdots t_{nn})w_1\cdots
w_n= \pi(t_{11})\cdots \pi(t_{nn})w_1\cdots w_n=w_1\cdots w_n.
\end{equation*}
Also $\pi(\det)w_1\cdots w_n=w_1\cdots w_n\otimes
\pair{\det}{\det}T_{-\underline 1}$, so
$\pair{\det}{\det}=T_{-2\cdot\underline 1}$.
\end{proof}

\begin{lemma}\label{lem:pairingGLn}
Define the pairing $\pair{\cdot}{\cdot}:\F_R(GL(n))^\cop\times
\F_R(GL(n))\to \Dh$ on the generators of $\F_R(GL(n))$ by
\eqref{eq:explicitpairingtij} and
\begin{equation}\label{eq:pairingdet-1}
\pair{\det^{-1}}{t_{ij}}=\delta_{ij}T_{\omega(\hat\imath)}, \quad
\pair{t_{ij}}{\det^{-1}}=\delta_{ij}T_{\omega(\hat\imath)}, \quad
\pair{\det^{-1}}{\det^{-1}}=T_{2\cdot \underline 1}.
\end{equation}
Then $\F_R(GL(n))^\cop$ and $\F_R(GL(n))$ are paired as
$\h$-bialgebroids.
\end{lemma}

\begin{proof}
For the pairing $\pair{\cdot}{\cdot}:\F_R(M(n))^\cop\times
\F_R(M(n))\to \Dh$ the statement follows from \cite{r}. Since
\begin{equation*}
\begin{split}
\delta_{ij}T_{-\omega(i)}&=\varepsilon(t_{ij})=\pair{t_{ij}}{1}=
\pair{t_{ij}}{\det\det^{-1}}
 =\sum_k\pair{t_{kj}}{\det}T_{\omega(k)}
\pair{t_{ik}}{\det^{-1}}=T_{-\underline 1}\pair{t_{ij}}{\det^{-1}},
\end{split}
\end{equation*}
the pairing is also well-defined for $\F_R(GL(n))$.
\end{proof}

We want to extend Lemma \ref{lem:pairingGLn} and show that the
pairing exists on the level of $\h$-Hopf algebroids. For this we
need to calculate pairings with dynamical quantum minor determinants
because of Proposition \ref{prop:dynMnantipode}, the proof of the
following lemma follows the same strategy as the proof of Lemma
\ref{lem:pairingtijdet}.

\begin{lemma}\label{lem:pairingtxi}
For $i\neq j$ we have
\begin{equation*}
\begin{split}
  &\pair{t_{ii}}{\qmindet{\imath}{\imath}}= q^{-1+1/n}\prod_{k<i}g(\lambda_k-\lambda_i) T_{-\underline
  1},\qquad
  \pair{t_{ii}}{\qmindet{\jmath}{\jmath}}= q^{1/n} T_{-\omega(\hat \jmath)-\omega(i)}, \\
  &\pair{t_{ij}}{\qmindet{\imath}{\jmath}}= q^{-1+1/n}h_0(\lambda_j-\lambda_i)
\prod_{k<j, k\neq i}-h(\lambda_j-\lambda_k)\prod_{k<i,k\neq j}
-h(\lambda_k-\lambda_i)T_{-\underline 1},
\end{split}
\end{equation*}
and
\begin{equation*}
\begin{split}
 &\pair{\qmindet{\imath}{\imath}}{t_{ii}}= q^{-1+1/n}\prod_{k>i}g(\lambda_i-\lambda_k) T_{-\underline
 1},\qquad
 \pair{\qmindet{\jmath}{\jmath}}{t_{ii}}=q^{1/n} T_{-\omega(\hat \jmath)-\omega(i)},\\
&\pair{\qmindet{\imath}{\jmath}}{t_{ij}}=q^{-1+1/n}
h_0(\lambda_i-\lambda_j) \prod_{k>j,k\neq i} -h(\lambda_k-\lambda_j)
\prod_{k>i,k\neq j}-h(\lambda_i-\lambda_k) T_{-\underline 1},
\end{split}
\end{equation*}
All other pairings between generators $t_{ij}$ and dynamical quantum
minor determinants $\qmindet{k}{l}$ are zero.
\end{lemma}

\begin{prop}
$\F_R(GL(n))^{\cop}$ and $\F_R(GL(n))$ are paired as $\h$-Hopf
algebroids.
\end{prop}

\begin{proof}
In Lemma \ref{lem:pairingGLn} we proved that $\F_R(GL(n))^\cop$ and
$\F_R(GL(n))$ are paired as $\h$-bialgebroids. So it remains to
check \eqref{eq:pairingforhHopfalgebras} on generators. The only non
trivial cases are $(X,a)= (t_{ii},t_{ii})$, $(t_{jj},t_{ii})$,
$(t_{ji},t_{ij})$, $(t_{ii}, \det^{-1})$ and $(\det^{-1},t_{ii})$
for $i\neq j$. From Example \ref{ex:oppcoop} we know $S^\cop=S^{-1}$
and $\mu_{l}^{\cop}=\mu_{r}$, $\mu_{r}^{\cop}=\mu_{l}$, so
\begin{equation*}
S^\cop(t_{ij})=\det^{-1}\xi^{\hat\jmath}_{\hat\imath}
\frac{\mu_r^\cop(\sign(j;\hat\jmath))}{\mu_l^\cop(\sign(i;\hat\imath))}.
\end{equation*}
Now we can check the relations by direct computations, using Lemma
\ref{lem:pairingtxi}. We show the third relation in detail; the
other cases can be done analogously. Using Lemma
\ref{lem:pairingtxi},
\begin{equation*}
\begin{split}
\pair{S^{\cop}(t_{ji})}{t_{ij}}& =
\pair{\det^{-1}\xi^{\hat\imath}_{\hat\jmath}
\frac{\mu_r^\cop(\sign(\{i\};\hat\imath))}{\mu_l^\cop(\sign(\{j\};\hat\jmath))}}{t_{ij}}\\
 &=\sign(\{j\};\hat\jmath)(\ulambda-\omega(\hat\jmath))^{-1}\pair{\det^{-1}}{t_{ii}}
 T_{\omega(i)}
 \pair{\xi^{\hat\imath}_{\hat\jmath}}{t_{ij}}\sign(\{i\};\hat\imath)(\ulambda)\\
  &= q^{-1+1/n}h_0(\lambda_i-\lambda_j) h(\lambda_i-\lambda_j)\prod_{m\neq
i,j} -h(\lambda_m-\lambda_j) \prod_{m\neq i,j}
-h(\lambda_i-\lambda_m),
\end{split}
\end{equation*}
and
\begin{equation*}
\begin{split}
 S(\pair{t_{ji}}{S(t_{ij})})&=S(\pair{t_{ij}}{\det^{-1}
\frac{\mu_l(\sign(\hat\jmath;\{j\}))}{\mu_r(\sign(\hat\imath;\{i\}))}\xi^{\hat\jmath}_{\hat\imath}})\\
 &= \sign(\hat\jmath;\{j\})(\ulambda)\pair{t_{ii}}{\det^{-1}} T_{\omega(i)}
\pair{t_{ji}}{\xi^{\hat\jmath}_{\hat\imath}}\sign(\hat\imath;\{i\})(\ulambda+\omega(\hat\imath))^{-1}\\
 &=q^{-1+1/n}h_0(\lambda_i-\lambda_j) h(\lambda_i-\lambda_j)\prod_{m\neq
i,j} -h(\lambda_m-\lambda_j) \prod_{m\neq i,j}
-h(\lambda_i-\lambda_m),
\end{split}
\end{equation*}
so $\pair{S^{\cop}(t_{ji})}{t_{ij}}=S(\pair{t_{ji}}{S(t_{ij})})$.
\end{proof}

\subsection{Compatible $\ast$-structures for the pairing}

If we equip $\F_R(GL(n))^\cop$ and $\F_R(GL(n))$ with the
$*$-operator defined in Lemma \ref{lem:dynGLn*structure}, they are
not paired as $\h$-Hopf $*$-algebroids. But since the $*$-operator
is not unique it is possible that there exists another $*$-operator
which gives paired $\h$-Hopf $*$-algebroids.

\begin{lemma}\label{lem:dynGLnpairedHopf}
The $\h$-Hopf algebroid $\F_R(GL(n))$ has a $*$-operator, denoted by
$\da$, defined on the generators by $\mu_l(f)^\da=\mu_l(\bar f)$,
$\mu_r(f)^\da=\mu_r(\bar f)$ and
\begin{equation}
    t_{ij}^\da=\frac{\mu_l(s_i)}{\mu_r(s_j)}
\xi^{\ihat}_{\jhat} \det^{-1},   \qquad (\det^{-1})^\da=\det,
\end{equation}
where $s_i(\ulambda)=q^{2(1/n\sum_{k=1}^n \lambda_k-\lambda_i)}$ and
extended as an $\C$-antilinear algebra anti-homo\-mor\-phism.
\end{lemma}

\begin{proof}
The proof follows the lines of the proof of Lemma
\ref{lem:dynGLn*structure}. On dynamical quantum minor determinants
we have
\begin{equation*}
(\xi^I_J)^\da= \frac{\mu_l(s_I)}{\mu_r(s_J)} \xi^{I^c}_{J^c}
\det^{-1},
\end{equation*}
where $s_I(\lambda)= q^{2(\# I/n\sum_{k=1}^n \lambda_k-\sum_{i\in I}
\lambda_i)}$. This follows  using $s_{I\setminus\{i\}}s_{\{i\}}=s_I$.
{}From the claim it follows that $\da$ is an involution. Indeed, since
$\mu_{l/r}(s_I)\: \det=\det\:\mu_{l/r}(s_I)$  and
$s_{\{1,\ldots,n\}}=1$ we have
\begin{equation*}
\begin{split}
(t_{ij}^\da)^\da &= \left(\frac{\mu_l(s_i)}{\mu_r(s_j)}
\xi^{\hat\imath}_{\hat\jmath} \det^{-1}\right)^\da= \det
\frac{\mu_l(s_{\hat\imath})}{\mu_r(s_{\hat\jmath})}t_{ij} \det^{-1}
\frac{\mu_l(s_i)}{\mu_r(s_j)}=t_{ij}.
\end{split}
\end{equation*}
Since the $*$-operator $*$ preserves the commutation relations and
$(t_{ij})^\da=\frac{\mu_l(s_i)}{\mu_r(s_j)} (t_{ij})^*$ it follows
directly from $R_{ab}^{xy}=0$ if $\omega(x)+\omega(y)\neq
\omega(a)+\omega(b)$ that $\da$ preserves the RLL-relations. By
direct computations we can check $\varepsilon\circ\da=
*^{\Dh}\circ\varepsilon$ and $(\da\otimes
\da)\circ\Delta=\Delta\circ\da$  on the generators and so on
$\F_R(GL(n))$.
\end{proof}

\begin{thm}\label{thm:dualGLn}
$(\F_R(GL(n))^{\cop},\da)$ and $(\F_R(GL(n)),*)$ are paired as
$\h$-Hopf $*$-alge\-broids.
\end{thm}

\begin{proof}
{}From Lemma \ref{lem:dynGLnpairedHopf} it follows that it remains to
prove \eqref{eq:pairingXstara} for generators. We have to check this
relation for five non-trivial cases: $(X,a)=(t_{ii},t_{ii})$,
$(t_{ii},t_{jj})$, $(t_{ij},t_{ij})$, $(\det^{-1},t_{ii})$ and
$(t_{ii},\det^{-1})$ for $i\neq j$. We give the proof of the second
case, which is the most involved one, in detail; the others can be
proved analogously. Since $T_{\omega(\hat\imath)}s_i(\ulambda)=
s_i(\ulambda)q^{2(1-1/n)}$, $T_{\omega(\hat\jmath)}s_i(\ulambda)=
s_i(\ulambda)q^{-2/n}$, for $i\neq j$, we have
\begin{equation*}
\begin{split}
\pair{t_{ii}^{\da}}{t_{jj}}&=s_i(\ulambda)^{-1}
\pair{\xi^{\hat\imath}_{\hat\imath}}{t_{jj}} T_{\omega(j)}
\pair{\det^{-1}}{t_{jj}} s_i(\ulambda)q^{2(1-1/n)}\\
 &= s_i(\ulambda)^{-1} q^{1/n}T_{-\omega(\hat\imath)+ \omega(\hat\jmath)}
s_i(\ulambda)q^{2(1-1/n)}=q^{-1/n} T_{\omega(i)-\omega(j)}.
\end{split}
\end{equation*}
For $i<j$ we also have
\begin{equation*}
\begin{split}
\pair{t_{ii}}{S(t_{jj})^*}& = \pair{t_{ii}}{t_{jj}
\frac{\mu_l(\sign(\hat\jmath;\{j\}))}{\mu_r(\sign(\hat\jmath;\{j\}))}}\\
 & =\sign(\hat\jmath;\{j\})(\ulambda-\omega(j))q^{-1/n}
T_{-\omega(i)-\omega(j)} \sign(\hat\jmath;\{j\})(\ulambda)^{-1} \\
 & =q^{-1/n}\frac{\sign(\hat\jmath;\{j\})(\ulambda-\omega(j))}
{\sign(\hat\jmath;\{j\})(\ulambda-\omega(i)-\omega(j))}T_{-\omega(i)-\omega(j)}
=q^{-1/n}T_{-\omega(i)-\omega(j)},
\end{split}
\end{equation*}
since $\sign(\hat\jmath;\{j\})$ is independent of $\lambda_i$ if
$i<j$. For $i>j$ we get
\begin{equation*}
\begin{split}
\pair{t_{ii}}{S(t_{jj})^*}& =
\sign(\hat\jmath;\{j\})(\ulambda-\omega(j))q^{-1/n}g(\lambda_j-\lambda_i)
T_{-\omega(i)-\omega(j)} \sign(\hat\jmath;\{j\})(\ulambda)^{-1}\\
 &=q^{-1/n}g(\lambda_j-\lambda_i)\frac{\sign(\hat\jmath;\{j\})(\ulambda-\omega(j))}
{\sign(\hat\jmath;\{j\})(\ulambda-\omega(i)-\omega(j))}T_{-\omega(i)-\omega(j)}\\
 &=q^{-1/n}g(\lambda_i-\lambda_j)\frac{h(\lambda_i-\lambda_j+1)}{h(\lambda_i-\lambda_j)}
T_{-\omega(i)-\omega(j)}=q^{-1/n}T_{-\omega(i)-\omega(j)},
\end{split}
\end{equation*}
where we use \eqref{eq:relationshg} in the last equality. So
$T_{-\omega(j)}\pair{t_{ii}}{S(t_{jj})^*}^*T_{-\omega(j)}= q^{-1/n}
T_{\omega(i)-\omega(j)}$ which proves the second case.
\end{proof}

\begin{remark}
Instead of the relation \eqref{eq:pairingXstara} we can also require
the pairing and $*$-operator to satisfy a similar relation where $*$
and $S$ are interchanged in the right hand side, see \cite{KvN}.
Also with that relation, the cobraiding \eqref{eq:cobrdynMn} on the
dynamical $GL(n)$ quantum group is not a pairing on the level of
$\h$-Hopf algebroids with the same $*$-operator $*$ on
$\F_R(GL(n))^\cop$ and $\F_R(GL(n))$.
\end{remark}



\begin{thebibliography}{999}

\bibitem{Aske} Askey, R., \emph{Orthogonal Polynomials and Special Functions}, Reg. Conf. Series in Appl. Math. {\bf 21}, SIAM, 1975.

\bibitem{BabeBB} Babelon, O., Bernard, D. and Billey, E.,
\emph{A quasi-Hopf algebra interpretation of quantum $3$-$j$ and
$6$-$j$ symbols and difference equations},  Phys. Lett. B
{\bf 375}  (1996),  89--97.

\bibitem{Bergman} Bergman, G.M., \emph{The diamond lemma for ring
theory}, Adv. in Math. {\bf 29} (1978), 178--218.

\bibitem{BrzeM} Brzezi\'nski, T. and Militaru, G.,
\emph{Bialgebroids, $\times_A$-bialgebras and duality},
J. Algebra {\bf 251} (2002), 279--294.

\bibitem{BuffR} Buffenoir, E. and Roche, Ph., 
\emph{An infinite product formula for $U_q(\rm{sl}(2))$
dynamical coboundary element},
J. Phys. A  {\bf 37}  (2004), 337--346.

\bibitem{Dijk} Dijkhuizen, M.S.,  \emph{Some remarks on the construction of quantum symmetric spaces},  Acta Appl. Math.  {\bf 44}  (1996),  59--80.

\bibitem{DijkK} Dijkhuizen, M.S., and  Koornwinder, T.H., 
\emph{CQG algebras: a direct algebraic approach to compact quantum groups}, 
Lett. Math. Phys. {\bf 32} (1994), 315--330.

\bibitem{DijkN} Dijkhuizen, M.S. and Noumi, M., \emph{A family of quantum projective spaces and related $q$-hypergeometric orthogonal polynomials}, Trans. Amer. Math. Soc.  {\bf 350}  (1998), 3269--3296.

\bibitem{EtingofSchiffmann} Etingof, P. and Schiffmann, O.,
\emph{Lectures on the dynamical Yang-Baxter equations}, pp. 89--129
in ``Quantum Groups and Lie Theory'', London Math. Soc. Lecture Note
Ser. Vol. 290, Cambridge Univ. Press, Cambridge, 2001, (extended
version at  {\tt math.QA/9908064}).

\bibitem{EtingofVarchenko}Etingof, P. and Varchenko, A., \emph{Solutions of the quantum dynamical
{Y}ang-{B}axter equation and dynamical quantum groups}, Comm. Math.
Phys. {\bf 196} (1998), 591-640.

\bibitem{EtinV2} Etingof, P. and Varchenko, A.,
\emph{Exchange dynamical quantum groups}, Comm. Math. Phys., {\bf
205} (1999), 19-52.

\bibitem{FelderVarchenko}Felder, G. and Varchenko, A., \emph{On representations of the elliptic quantum group ${E}\sb
{\tau,\eta}({\rm sl}\sb 2)$}, Comm. Math. Phys. {\bf 181} (1996),
741-761.

\bibitem{Flor}  Floris, P.G.A., \emph{Addition formula for $q$-disk polynomials}, 
Comp. Math. {\bf 108} (1997), 123--149.

\bibitem{Kass} Kassel, C., \emph{Quantum Groups}, GTM 155, Springer, 1995.
 
\bibitem{KoelAAM} Koelink, H.T., \emph{Askey-Wilson polynomials and the quantum ${\rm SU}(2)$ group: survey and applications},  Acta Appl. Math.  
{\bf 44}  (1996), 295--352. 

\bibitem{KoelFI} Koelink, E.,  \emph{Addition formulas for $q$-special functions},  109--129 in ``Special Functions, $q$-series and Related Topics'', 
(eds. M.E.H. Ismail, D.R. Masson and M. Rahman), Fields Inst. Commun. {\bf 14}, AMS, 1997.

\bibitem{KvN} Koelink, E. and Norden, Y.~van, \emph{Pairings and actions for
dynamical quantum groups}, {\tt math.QA/0412205}.

\bibitem{KvNR} Koelink, E., Norden, Y. van, and Rosengren, H.,
\emph{Elliptic $U(2)$ quantum group and elliptic
hypergeometric series},
Comm. Math. Phys.  {\bf 245}  (2004), 519--537.
(extended version at  {\tt math.QA/0304189})

\bibitem{KoelinkRosengren}Koelink, E. and Rosengren, H.,\emph{Harmonic analysis on the {S}{U}$(2)$ dynamical quantum
group}, Acta Appl. Math. {\bf 69} (2001), 163-220.

\bibitem{KoorIM} Koornwinder, T. H.,\emph{The addition formula for Jacobi polynomials. I. Summary of results}, Indag. Math.  {\bf 34}  (1972),  188--191. 

\bibitem{Koor} Koornwinder, T.H., \emph{Askey-Wilson polynomials as zonal spherical functions on the ${\rm SU}(2)$ quantum group},  SIAM J. Math. Anal.  {\bf 24}
(1993), 795--813.

\bibitem{Letz} Letzter, G., \emph{Quantum zonal spherical functions and
Macdonald polynomials}, Adv. Math. {\bf 189} (2004),
88--147.

\bibitem{Lu} Lu, J.-H., \emph{Hopf algebroids and quantum groupoids},
Internat. J. Math. {\bf 7} (1996), 47--70.

\bibitem{Macdonald} Macdonald, I. G., \emph{Symmetric functions and {H}all polynomials},
Oxford Mathematical Monographs, The Clarendon Press Oxford
University Press, New York, 1995.

\bibitem{McCoR}  McConnell, J.C. and  Robson, J.C., 
\emph{Noncommutative Noetherian Rings}, John Wiley, 1987.

\bibitem{vNord} Norden, Y. van, \emph{Dynamical Quantum Groups:
Duality and Special Functions}, thesis, TU Delft, 2005.

\bibitem{Noum} Noumi, M., \emph{ Macdonald's symmetric polynomials as zonal spherical functions on some quantum homogeneous spaces}
Adv. Math. {\bf 123} (1996),  16--77.

\bibitem{NoumM} Noumi, M. and Mimachi, K., \emph{Askey-Wilson polynomials and the quantum group ${\rm SU}_q(2)$},  Proc. Japan Acad. Ser. A Math. Sci.  {\bf 66}  (1990),  146--149

\bibitem{NoumiYM}Noumi, M., Yamada, H. and Mimachi, K., \emph{Finite-dimensional 
representations of the quantum group {${\rm GL}\sb q(n;{\bf C})$}
and the zonal spherical functions on {${\rm U}\sb
q(n-1)\backslash{\rm U}\sb q(n)$}}, Japan. J. Math. (N.S.) {\bf 19}
(1993), {31-80}.

\bibitem{RoseCM} Rosengren, H., \emph{A new quantum algebraic 
interpretation of the Askey-Wilson polynomials}, Contemp. Math. {\bf 254}  (2000), 371--394.

\bibitem{r}Rosengren, H., \emph{Duality and self-duality for dynamical quantum
groups}, Algebr. Represent. Theory, {\bf 7} (2004), no. 4, 363--393.

\bibitem{Stok} Stokman, J.V.,
\emph{Vertex-IRF transformations, dynamical quantum groups
and harmonic analysis},
Indag. Math. N.S. {\bf 14} (2003), 545--570.

\bibitem{StokmanOb}Oblomkov, A.A. and Stokman, J.V., \emph{Vector valued spherical func\-tions and
Macdonald-\-Koornwinder polynomials}, 
Comp. Math., to appear,  {\tt math.QA/0311512}.

\bibitem{Take} Takeuchi, M., \emph{Groups of algebras over $A\otimes\bar A$},
J. Math. Soc. Japan {\bf 29} (1977), 459--492.

\bibitem{Xu} Xu, P., \emph{Quantum groupoids}, Comm. Math. Phys. {\bf 216}
(2001), 539--581.

\end{thebibliography}
\end{document}